\documentclass{article}
\usepackage[utf8]{inputenc}
\usepackage[T1]{fontenc} 
\usepackage[linesnumbered,ruled,vlined]{algorithm2e} 
\usepackage{amsmath}
\usepackage{amsfonts}
\usepackage{amssymb}
\usepackage{amsthm}
\usepackage{amscd}
\usepackage{color}
\usepackage{esint}
\usepackage{enumerate}
\usepackage{dsfont}
\usepackage{psfrag}
\usepackage{stmaryrd}
\usepackage{hyperref}
\usepackage{graphicx}
\usepackage{graphics}
\usepackage{epstopdf}
\usepackage{epsfig}
\usepackage{setspace}
\usepackage{tikz}
\usetikzlibrary{positioning,arrows}
\usetikzlibrary{calc,arrows}
\usetikzlibrary{er,positioning,bayesnet}
\usepackage{pgfplots}
\usepackage{pgfplotstable}
\pgfplotsset{compat=1.11} 
\usepackage{bm}
\usepackage{booktabs} 
\usepackage{caption} 
\usepackage{subcaption} 
\usepackage[all]{nowidow}
\usepackage{comment} 
\usepackage{multicol}
\usepackage{mathabx}
\usepackage{mathrsfs}
\usepackage{mathtools}
\usepackage{fullpage}

\usepackage[inline]{enumitem} 

\newtheorem{theorem}{Theorem}

\newtheorem{remark}[theorem]{Remark}

\newcommand{\RR}{\mathbb{R}}

\newcommand{\dps}{\displaystyle}
\newcommand{\dsp}{\displaystyle}

\def\F{\mathcal{F}}
\def\C{\mathcal{C}}
\def\dt{{\Delta t}}
\def\ddt{{\delta t}}

\DeclareMathOperator{\cost}{Cost}

\def\dc{\delta_{\mathrm{conv}}}
\def\de{\delta_{\mathrm{expl}}}

\def\Nin{N_{\mathrm{init}}}
\def\Nfi{N_{\mathrm{final}}}

\def\Nad{N_{\mathrm{slab}}}

\def\qprev{q^{\mathrm{prev}}}
\def\qcur{q^{\mathrm{cur}}}
\def\pprev{p^{\mathrm{prev}}}
\def\pcur{p^{\mathrm{cur}}}

\usepackage{shadethm}

\definecolor{blue}{HTML}{1F77B4}
\definecolor{orange}{HTML}{FF7F0E}
\definecolor{green}{HTML}{2CA02C}

\begin{document}
\date{\today}
\author{Olga Gorynina$^{1,2}$, Frédéric Legoll$^{3,2}$, Tony Lelièvre$^{1,2}$ and Danny Perez$^4$\\
{\footnotesize $^1$ CERMICS, École des Ponts,}\\
{\footnotesize 6 et 8 avenue Blaise Pascal, 77455 Marne-La-Vall\'ee Cedex 2, France}\\
{\footnotesize $^2$ MATHERIALS project-team, Inria Paris,}\\
{\footnotesize 2 rue Simone Iff, CS 42112, 75589 Paris Cedex 12, France}\\
{\footnotesize $^3$ Navier, École des Ponts, Univ Gustave Eiffel, CNRS,}\\
{\footnotesize 6 et 8 avenue Blaise Pascal, 77455 Marne-La-Vall\'ee Cedex 2, France}\\
{\footnotesize $^4$ Theoretical Division T-1, Los Alamos National Laboratory,}\\
{\footnotesize Los Alamos, NM 87545, USA} \\
{\footnotesize \tt \{olga.gorynina,frederic.legoll,tony.lelievre\}@enpc.fr, danny\_perez@lanl.gov}
}
\title{Combining machine-learned and empirical force fields with the parareal algorithm: application to the diffusion of atomistic defects}

\maketitle

\begin{abstract}
We numerically investigate an adaptive version of the parareal algorithm in the context of molecular dynamics. This adaptive variant has been originally introduced in~\cite{papierUpanshu}. We focus here on test cases of physical interest where the dynamics of the system is modelled by the Langevin equation and is simulated using the molecular dynamics software LAMMPS. In this work, the parareal algorithm uses a family of machine-learning spectral neighbor analysis potentials (SNAP) as fine, reference, potentials and embedded-atom method potentials (EAM) as coarse potentials. We consider a self-interstitial atom in a tungsten lattice and compute the average residence time of the system in metastable states. Our numerical results demonstrate significant computational gains using the adaptive parareal algorithm in comparison to a sequential integration of the Langevin dynamics. We also identify a large regime of numerical parameters for which statistical accuracy is reached without being a consequence of trajectorial accuracy.
\end{abstract}


\section{Introduction}



This work is motivated by molecular dynamics (MD) simulations, where one often has to compute ensemble averages or dynamical quantities, which both involve averages over very long trajectories of stochastic dynamics (we refer e.g. to~\cite{LelievreRoussetStoltz10} for a general, mathematically oriented presentation of the MD context). Reducing the computational cost of these long time simulations, or at least the wall-clock time it takes to obtain such long trajectories, is thus of great practical interest. As conventional spatial parallelization schemes based on domain decomposition allow for larger system sizes but not longer simulation times~\cite{Uberuaga2018}, one alternative way to speed up such computations is to design accelerated MD approaches based on the parallelization of the problem in the temporal domain~\cite{zamora2020accelerated}.

A popular parallel-in-time method for integrating ordinary differential equations is the parareal algorithm, which has been first introduced in~\cite{lions2001resolution}. Using a time domain decomposition, the algorithm aims at computing, in an iterative manner, an approximation of the exact solution of the dynamics (see Section~\ref{sec:parareal} for a detailed description). The whole time domain is divided into several subintervals. At each iteration, the parareal algorithm utilizes a coarse solver to quickly step through the time domain by computing relatively cheap approximate solutions for all time subintervals, and then simultaneously refines all of these approximate solutions using an accurate fine solver which is applied {\em in parallel} over each time subinterval. Since the fine propagator corrections (which are expensive to compute) are applied concurrently over the subintervals (and not in a sequential manner from the initial until the final time), reduction in the associated wall-clock time is possible. In contrast, the coarse propagator is applied sequentially over the complete time interval, but its cost is often negligible when compared with the cost of the fine propagator.

In the MD context, it is convenient to consider parareal algorithms where the coarse and the fine propagators integrate dynamics based on different potential energies, using the same time-step, the difference in cost stemming from the different complexity for evaluating the potential (indeed, most physically relevant coarse potentials would likely have very similar numerical-stability constraints as fine potentials). In this setting, an adaptive version of the parareal algorithm has recently been introduced in~\cite{papierUpanshu}. It is shown there that this adaptive version leads to significantly improved gains (in comparison to the standard version of parareal) on some toy examples. 

The main goal of this article is to investigate the performances of the adaptive parareal algorithm for realistic problems of physical interest. To that aim, we focus on the MD simulation of the diffusion of a self interstitial atom (SIA) in a body-centered cubic (BCC) tungsten lattice. The simulation is performed using the LAMMPS~\cite{LAMMPS} molecular dynamics software (Large-scale Atomic/Molecular Massively Parallel Simulator), a software which is very broadly used within the materials science community. To model the tungsten lattice, we consider several interatomic potentials in two families: machine-learned spectral neighbor analysis potentials (SNAP)~\cite{thompson2015spectral} and embedded-atom method potentials (EAM)~\cite{daw1984embedded}.

This article is organized as follows. In Section~\ref{sec:algo}, we briefly review the classical parareal algorithm and the adaptive version introduced in~\cite{papierUpanshu}. We also discuss the implementation of our method in LAMMPS, including some issues related to the fact that we use {\em as is} the time-stepping scheme implemented in LAMMPS (for the sake of having a non-intrusive implementation). This raises some unexpected difficulties (some calculation details related to these questions are postponed until Appendix~\ref{sec:app}). We next present our SIA simulations in Section~\ref{sec:num}. We describe the MD settings we have chosen and demonstrate the accuracy (both from a trajectorial and a statistical viewpoint) of the results obtained using the adaptive parareal algorithm in comparison to reference results obtained using a standard sequential integration of the dynamics. We also discuss in Section~\ref{sec:num} the significant computational gains obtained when using the adaptive parareal algorithm (on our test cases, this gain varies between 3 and 20 depending on the choice of the coarse propagator and of the time-step).

\section{Algorithm} \label{sec:algo}

The present work focuses on the parallelization in time of the Langevin dynamics
\begin{equation} \label{eq:langevin}
    \begin{cases}
    dq(t) = p(t) \, dt, \\
    dp(t) = - \nabla V(q(t)) \, dt - \gamma \, p(t) \, dt + \sqrt{2 \gamma \beta^{-1}} \, d W(t),
    \end{cases}
\end{equation}
with initial condition $(q_0,p_0)$. Here $(q(t),p(t)) \in \RR^{3m} \times \RR^{3m}$ are the positions and momenta at time $t \in [0,T]$ of $m$ particles in the three-dimensional space, $V : \RR^{3m} \rightarrow \RR$ is the potential energy of the system, $\gamma > 0$ is the friction coefficient, $\beta$ is (up to a multiplicative constant) the inverse temperature and $W(t)$ is a standard Brownian motion in dimension $3m$. In~\eqref{eq:langevin}, we have set the mass of each particle to unity for simplicity, but the generalization to non-identity mass matrix is straightforward (and the numerical tests discussed in Section~\ref{sec:num} are actually performed with the physically relevant values of the particle masses).

Let us introduce a uniform grid on the time interval $[0,T]$ with $N$ time-steps of length $\dt = T/N$: 
$$
0 = t_0 <t_1 = \dt <\dots < N\dt = t_N = T.
$$ 
Because it will be useful in the context of the parareal algorithm, we consider that each time-step $\dt$ may be itself subdivided into $L$ time-steps of length $\ddt = \dt/L$ for some $L \geq 1$. 

The integrator for the Langevin dynamics~\eqref{eq:langevin} implemented in LAMMPS is a modification of the BBK scheme~\cite{brunger1984stochastic} with an effective force which takes into account the damping term and the fluctuation term associated to the white noise, in addition to the physical force $-\nabla V$. Let us consider $L+1$ independent standard Gaussian random variables, denoted $G_0$, $G_1$, \dots, $G_L$. The first iterate of the scheme has a particular expression:
\begin{equation} \label{eq:scheme_debut}
    \begin{cases}
        p_{1/2} = p_0 - \cfrac{\ddt}{2} \, \nabla V(q_0) - \cfrac{\ddt}{2} \, \gamma \, p_0 + \cfrac{1}{2} \, \sqrt{2 \gamma \beta^{-1} \ddt} \, G_0,\\
        q_1 = q_0 + \ddt \, p_{1/2}, \\
        p_1 = p_{1/2} - \cfrac{\ddt}{2} \, \nabla V(q_1) - \cfrac{\ddt}{2} \, \gamma \, p_{1/2} + \cfrac{1}{2} \, \sqrt{2 \gamma \beta^{-1} \ddt} \, G_1,
    \end{cases}
\end{equation}
where we recall that $(q_0, p_0)$ is the initial condition. The subsequent iterates (for $1 \leq \ell \leq L-1$) are given by
\begin{equation} \label{eq:scheme}
    \begin{cases}
        p_{\ell+1/2} = p_\ell - \cfrac{\ddt}{2} \, \nabla V(q_\ell) - \cfrac{\ddt}{2} \, \gamma \, p_{\ell-1/2} + \cfrac{1}{2} \, \sqrt{2 \gamma \beta^{-1} \ddt} \, G_\ell,\\
        q_{\ell+1} = q_\ell + \ddt \, p_{\ell+1/2}, \\
        p_{\ell+1} = p_{\ell+1/2} - \cfrac{\ddt}{2} \, \nabla V(q_{\ell+1}) - \cfrac{\ddt}{2} \, \gamma \, p_{\ell+1/2} + \cfrac{1}{2} \, \sqrt{2 \gamma \beta^{-1} \ddt} \, G_{\ell+1}.
    \end{cases}
\end{equation}
The first iterate~\eqref{eq:scheme_debut} differs from the next ones~\eqref{eq:scheme} in two ways. First, if we were to set $\ell=0$ in the first line of~\eqref{eq:scheme}, we would need to know $p_{-1/2}$ to compute $p_{1/2}$. However, $p_{-1/2}$ is of course not defined. It is replaced in the first line of~\eqref{eq:scheme_debut} by $p_0$. Second, and more importantly, the random variable $G_0$ is only used in the first iteration, whereas all the other random variables $\{ G_\ell \}_{1 \leq \ell \leq L-1}$ are used twice ($G_\ell$ is used in the last step of iterate $\ell-1$ and in the first step of iterate $\ell$). Stated otherwise, the first iterate~\eqref{eq:scheme_debut} uses {\em two} independent random variables, whereas each subsequent iterate only uses {\em one} additional independent random variable. This difference has far-reaching consequences, e.g. on the equilibrium kinetic temperature simulated by the numerical scheme, as will be discussed in Section~\ref{section:eff_temp}.

\begin{remark} \label{rem:pourquoi_implem}
The modification~\eqref{eq:scheme_debut}--\eqref{eq:scheme} of the BBK scheme is particularly well-suited to adapt to the Langevin equation an implementation of the standard Verlet algorithm used to integrate deterministic Hamiltonian dynamics. Indeed, as in a Verlet scheme, for $\ell \ge 1$, the total force (sum of force field, damping and fluctuation) to compute $p_{\ell+1/2}$ from $p_\ell$ in the first step of iterate $\ell$ is exactly the same as the total force to compute $p_\ell$ from $p_{\ell-1/2}$ in the last step of iterate $\ell-1$. The total force hence only needs to be computed once per time-step. 
\end{remark}

\subsection{Parareal method} \label{sec:parareal}

To integrate~\eqref{eq:langevin} over a time interval of length $\dt$, we consider two propagators: $\F_{\dt}$ is a fine, expensive propagator which accurately approximates the exact solution of~\eqref{eq:langevin}, and $\C_{\dt}$ is a coarse, less expensive propagator which is also less accurate. In the MD context, $\F_{\dt}$ and $\C_{\dt}$ are often integrators with a fixed discretization scheme and with the same time-step, but run on different potential energy landscapes. In what follows, and similarly to~\cite{papierUpanshu}, the fine propagator $\F_{\dt}$ performs $L$ iterations of~\eqref{eq:scheme} (with the time-step $\ddt = \dt/L$), where $V \equiv V_{\F}$ is the reference, accurate potential energy. The coarse propagator $\C_{\dt}$ also performs $L$ iterations of~\eqref{eq:scheme} (with the same time-step $\ddt$), but where $V \equiv V_{\C}$ is now a coarse, approximate potential energy. As explained above, for both propagators $\F_{\dt}$ and $\C_{\dt}$, the first of these $L$ iterates is given in the form of~\eqref{eq:scheme_debut} instead of~\eqref{eq:scheme}. We emphasize that, to ensure convergence of the parareal algorithm, both propagators $\F_{\dt}$ and $\C_{\dt}$ should use the {\em same} random variables. We denote by $\{ G_{\ell,n} \}_{0 \leq \ell \leq L, \, 1 \leq n \leq N}$ the random variables used to propagate the system from the initial time to the final time $N \, \dt$.

The classical parareal method (as first introduced in~\cite{lions2001resolution} and reformulated in~\cite{bal-maday-02}) is an iterative, parallel-in-time algorithm, which computes the trajectory over $[0,T]$ by using domain (in time) decomposition. It proceeds as follows. Starting from the initial condition $(q_0,p_0)$, the algorithm first performs a sequential coarse propagation to compute $\{ (q^0_n, p^0_n) \}_{0\leq n\leq N}$: 
\begin{equation} \label{eq:parareal_0}
    (q^0_{n+1}, p^0_{n+1}) = \C_{\dt}(q^0_n, p^0_n), \qquad (q^0_0, p^0_0) = (q_0,p_0).
\end{equation}
Suppose now that we have at hand some numerical trajectory $\{ (q^{k-1}_n, p^{k-1}_n) \}_{0\leq n\leq N}$, obtained at the previous iteration $k-1$. The new parareal solution $\{ (q^k_n, p^k_n) \}_{0\leq n\leq N}$ is computed from the following scheme:
\begin{equation} \label{eq:parareal}
    (q^k_{n+1}, p^k_{n+1}) = \C_{\dt}(q^k_n, p^k_n) + \F_{\dt}(q^{k-1}_n, p^{k-1}_n) -\C_{\dt}(q^{k-1}_n, p^{k-1}_n), \qquad (q^k_0, p^k_0) = (q_0,p_0).
\end{equation}
We thus propagate the system $(q^{k-1}_n, p^{k-1}_n)$, in parallel over the time-windows $[n\dt,(n+1)\dt]$, according to both the coarse and the fine propagators, thereby obtaining $\F_{\dt}(q^{k-1}_n, p^{k-1}_n)$ and $\C_{\dt}(q^{k-1}_n, p^{k-1}_n)$. We next compute, still in parallel, the jumps $\F_{\dt}(q^{k-1}_n, p^{k-1}_n) -\C_{\dt}(q^{k-1}_n, p^{k-1}_n)$. We next perform a sequential propagation, from the initial condition and using the coarse propagator that we correct according to the precomputed jumps. We note that the fine solver is only used in the parallel part of the algorithm (the fine solver is applied in each interval $[n\dt,(n+1)\dt]$ independently of the other intervals), while the sequential part of the algorithm only calls the coarse propagator. The random numbers $\{ G_{\ell,n} \}_{0 \leq \ell \leq L, \, 1 \leq n \leq N}$ which are used are the same at each parareal iteration $k$.


The parareal algorithm has been successfully applied to many problems. We refer to~\cite{gander2007analysis} for a reformulation of the algorithm in a more general setting that relates the parareal strategy to earlier time-parallel algorithms (see also~\cite{gander-lunet,gander_50years}). Several variants of the algorithm have been proposed for specific applications: multiscale-in-time problems (see e.g.~\cite{BBK,maday41parareal}, \cite{LegollLelievreSamaey13} and~\cite{engblom2009parallel}), Hamiltonian ODEs or PDEs~\cite{DaiBrisLegollMaday13,dai_maday}, stochastic differential equations~\cite{bal2003parallelization,sall2016,LegollLelievreSamaey20}, reservoir simulations~\cite{garrido2005}, applications in fluid and solid mechanics~\cite{farhat2003time}, to mention but a few. This work is based on the adaptive version of the algorithm that we introduced in~\cite{papierUpanshu}. We also wish to mention~\cite{MulaMaday20} for another work in the direction of designing adaptive variants of the parareal algorithm. 

\medskip

In the following, we denote by $\{ (q_n^{\rm ref}, p_n^{\rm ref}) \}_{n \geq 0}$ the reference trajectory obtained by a sequential use of the fine propagator $\F_\dt$ (and which is thus expensive to compute): $(q_{n+1}^{\rm ref}, p_{n+1}^{\rm ref}) = \F_\dt(q_n^{\rm ref}, p_n^{\rm ref})$. An important feature of the parareal method is that, on the interval $[0,N \dt]$, it is guaranteed to converge to this reference solution in at most $k = N$ iterations (see e.g.~\cite{bal-maday-02}): for any $k \geq N$, we have $(q^k_n, p^k_n) = (q_n^{\rm ref}, p_n^{\rm ref})$ for any $0 \leq n \leq N$. In practice, convergence is often reached much sooner, therefore providing computational gains (a noteworthy exception is the case of hyperbolic problems, where a larger number of iterations is often needed to reach convergence, as observed e.g. in~\cite{DaiBrisLegollMaday13,dai_maday}; we also refer to the analysis in~\cite{papierUpanshu}).

Let us introduce the relative error between the reference trajectory $\{ q_n^{\rm ref} \}_{0\leq n\leq N}$ and the parareal trajectory $\{ q^k_n \}_{0\leq n\leq N}$ at the iteration $k$:
$$
    E_{\rm ref}\left(q^k,N\right) = \frac{\sum_{n=1}^N |q_n^{\rm ref} - q^k_n|}{\sum_{n=1}^N |q_n^{\rm ref}|}.
$$
The error $E_{\rm ref}$ cannot be computed in practice, since we do not have access to the reference trajectory $\{ q_n^{\rm ref} \}_{n \geq 0}$. Therefore, in order to monitor the convergence of the parareal method along the iterations, we introduce the relative error between two consecutive parareal trajectories:
\begin{equation} \label{def:rel-Err}
    E\left(q^{k-1},q^k,N\right) = \frac{\sum_{n=1}^N |q^k_n - q^{k-1}_n|}{\sum_{n=1}^N |q^{k-1}_n|}.
\end{equation}
As mentioned above, we have $q^k_n = q_n^{\rm ref}$ for any $k \geq n$, and thus $E\left(q^{k-1},q^k,N\right) = 0$ for $k \geq N+1$.

In our numerical experiments discussed below, we proceed with the parareal iterations until the accuracy reaches the user-chosen threshold $\dc$, namely until $E\left(q^{k-1},q^k,N\right)<\dc$. We denote by $k_{\rm conv}$ the number of parareal iterations required to reach this accuracy. 

\begin{remark}
Note that other criteria of convergence can be envisioned. For example, one could imagine a criteria based on the statistical properties of the numerical solution $\{ (q^k_n, p^k_n) \}_{0 \leq n \leq N}$. We do not address this question in this work.
\end{remark}

The classical parareal algorithm, as described above, is presented as Algorithm~\ref{algo:classical}.

\medskip

\begin{algorithm}[h]
\SetAlgoLined
\SetKw{KwBy}{by}
\SetKwInput{Input}{Numerical parameters}\SetKwInOut{Output}{Output of the algorithm}
\Input{$N$, $\dc$}
Compute $\{(\qcur_n, \pcur_n)\}_{0\leq n\leq N}$: $(\qcur_0, \pcur_0) := (q_0,p_0)$ and $(\qcur_{n+1}, \pcur_{n+1}) := \C_\dt(\qcur_n, \pcur_n)$\;
Set $k := 0$ and $\delta := 2 \, \dc$\;
\While{$\delta \geq \dc$}
{	
	\BlankLine
	$k := k + 1$\;
	Define $\{(\qprev_n, \pprev_n)\}_{0\leq n\leq N}$ as $(\qprev_n, \pprev_n) := (\qcur_n, \pcur_n)$ for any $0\leq n\leq N$\; 
	Compute $\{J_n\}_{0\leq n\leq N-1}$ in parallel: $J_n := \F_\dt(\qprev_n, \pprev_n) - \C_\dt(\qprev_n, \pprev_n)$\;
    \For{$n\gets 0$ \KwTo $N-1$ \KwBy $1$}
	{
    	$(\qcur_{n+1}, \pcur_{n+1}) := \C_\dt(\qcur_n, \pcur_n) + J_n$\;
  	}
	Compute the relative error $\delta = E\left(\qprev,\qcur,N\right)$\;
}
\Output{$\{(\qcur_n, \pcur_n)\}_{0 \leq n \leq N}$ and $k_{\rm conv} := k$}
\caption{Parareal algorithm \label{algo:classical}}
\end{algorithm}

\medskip

Let us denote by $C_f$ (resp. $C_c$) the cost of a single evaluation of the fine integrator $\F_{\dt}$ (resp. coarse integrator $\C_{\dt}$). Assuming that the communication time is zero, the total cost of the parareal algorithm is $\dsp N \, C_c + k_{\rm conv} \big( (C_f+C_c) + N \, C_c \big)$. In contrast, the cost of a sequential propagation according to the fine propagator is $N \, C_f$. We are thus able to define the wall-clock gain of the parareal method as 
$$
\Gamma(\dc, N) = \frac{N \, C_f}{N \, C_c + k_{\rm conv} \big( (C_f+C_c) + N \, C_c \big)}.
$$
If we additionally assume that the cost of the coarse propagations is negligible in comparison with the cost of the fine propagator (i.e. $N \, C_c \ll C_f$), we get $\Gamma(\dc, N) = \Gamma^{\rm ideal}(\dc, N)$, where the ideal gain is defined by $\dsp \Gamma^{\rm ideal}(\dc, N) := \frac{N}{k_{\rm conv}}$. Note that the total CPU effort spent by the parareal algorithm, which is equal to $\dsp N \, C_c + k_{\rm conv} \big( N (C_f+C_c) + 2 \, N \, C_c \big)$, is of course larger than the total CPU effort spent by a sequential fine integration (which is equal to $N \, C_f$). The interest of the parareal algorithm is that the computational effort is performed {\em in parallel} over several processors, therefore potentially providing a gain in terms of wall-clock time. We also note that computational efficiency considerations can further increase the attractiveness of parareal. For example, MD simulations on small systems can be very inefficient on modern Graphical Processing Units (GPUs), given the extremely high level of data parallelism that can be accommodated by the hardware. In this case, parareal could be implemented to execute all the fine-grained steps simultaneously on the {\em same} GPU instead of requiring parallelization over multiple GPUs, which could significantly improve the efficiency of the calculation and lead to a decrease in the net computational effort required to carry out the simulation.

\subsection{Adaptive parareal method}

The work~\cite{papierUpanshu} has introduced an adaptive version of the parareal algorithm, motivated by the fact that, when applied to MD problems, the classical parareal algorithm suffers from various limitations (in particular, possible intermediate blow-up of the trajectory and lack of computational gain in the case of too long time horizons). The new approach introduced in~\cite{papierUpanshu} consists in adaptively dividing (on the basis of the relative error between two consecutive parareal trajectories) the computational domain $[0,N \, \dt]$ in several subdomains. For that, we have to revisit the definition of the error~\eqref{def:rel-Err} for an arbitrary time-slab $[\Nin \dt,\Nfi \dt]$, for some fixed $\Nfi \geq \Nin\geq 0$. We naturally extend~\eqref{def:rel-Err} as
$$
E\left(q^{k-1},q^k,\Nin,\Nfi\right) = \frac{\sum_{n=\Nin}^{\Nfi} |q^k_n - q^{k-1}_n|}{\sum_{n=\Nin}^{\Nfi} |q^{k-1}_n|}.
$$

The adaptive parareal method proceeds as follows (see Algorithm~\ref{algo:adap}). We fix two parameters $\dc > 0$ (convergence parameter) and $\de > 0$ (explosion threshold parameter, which satisfies $\de > \dc$) and start by running the classical parareal algorithm on the whole time-slab $[0,T]$. As in the classical parareal method, for every parareal iteration $k$, we check whether the trajectory has reached convergence on the whole time slab, i.e. whether
$$
E\left(q^{k-1}, q^k,0 , N\right) < \dc.
$$
If this is the case then we stop the iterations. If this is not the case, then
\begin{itemize}
    \item either $E\left(q^{k-1}, q^k,0 , N\right) \leq \de$, and we then proceed with the next parareal iteration over the whole time-slab $[0,T]$.
    \item or $E\left(q^{k-1},q^k, 0 , N\right) > \de$, which means that the relative error at the parareal iteration $k$ is too large. We then give up on trying to reach convergence on the whole time-slab $[0,T]$ and decide to shorten it. In practice, we look for the smallest $n$ such that $E\left(q^{k-1},q^k, 0 , n\right) > \de$, denote it $n^{\rm cur}$ and shorten the original time-slab $[0,N\dt]$ to $[0,n^{\rm cur} \, \dt]$. 
    
    We then proceed with the next parareal iterations on this new time-slab, that will possibly be further shortened, until the relative error $E$, on the shortened time-slab $[0,n^{\rm cur} \, \dt]$, becomes smaller than $\dc$. We have then reached convergence on the current time-slab and proceed with the subsequent part of the time range. We thus define the new (tentative) time-slab as $[n^{\rm cur} \, \dt, N \, \dt]$ and start again the adaptive parareal algorithm. This procedure is repeated until the final time $T$ is reached.
\end{itemize}
The adaptive parareal method, as described above, is presented as Algorithm~\ref{algo:adap}. We denote there by $\Nad$ the number of time-slabs in which the whole time range $[0,N \, \dt]$ is eventually divided: $\dps [0,N \, \dt] = \bigcup_{1 \leq i \leq \Nad} [\Nin^i \, \dt, \Nfi^i \, \dt]$ with $\Nin^1=0$, $\Nfi^i = \Nin^{i+1}$ and $\Nfi^{\Nad} = N$. For any $1 \leq i \leq \Nad$, we denote by $k^i_{\rm conv}$ the number of parareal iterations required to reach convergence on $[\Nin^i \, \dt, \Nfi^i \, \dt]$. 

\medskip

\begin{algorithm}[h]
\SetAlgoLined
\SetKw{KwBy}{by}
\SetKwInput{Input}{Numerical parameters}\SetKwInOut{Output}{Output of the algorithm}
\Input{$N$, $\dc$, $\de$}
Set $\Nin:=0$, $\Nfi := 0$, and $\delta := (\dc+\de)/2$\;
Set $(\qcur_0, \pcur_0) := (q_0,p_0)$\;
Set $\Nad := 0$\;
\While{$\Nfi<N$}
{	
	\BlankLine
    \If{$\delta < \de$}
	    {
            $\Nin := \Nfi$\tcp*[l]{initialization for a new time-slab}
            $\Nfi := N$\;
            Compute $\{(\qcur_n, \pcur_n)\}_{\Nin\leq n\leq \Nfi}$: $\forall \Nin \leq n \leq \Nfi-1$, $(\qcur_{n+1}, \pcur_{n+1}) := \C_\dt(\qcur_n, \pcur_n)$\; 
            $\Nad := \Nad + 1$\;
            $k^{\Nad}_{\rm conv} := 0$\;
        }
    \BlankLine    
    Set $\delta := (\dc+\de)/2$\;
	\While{$\delta\in [\dc,\de]$}
	{	
		\BlankLine
		Define $\{(\qprev_n, \pprev_n)\}_{\Nin\leq n \leq \Nfi}$ as $(\qprev_n, \pprev_n) := (\qcur_n, \pcur_n)$ for all $\Nin \leq n \leq \Nfi$\; 
		Compute $\{J_n\}_{\Nin\leq n \leq \Nfi-1}$ (in parallel): $J_n := \F_\dt(\qprev_n, \pprev_n) - \C_\dt(\qprev_n, \pprev_n)$\;
        $k^{\Nad}_{\rm conv} := k^{\Nad}_{\rm conv} + 1$\;
		\For{$n\gets \Nin$ \KwTo $\Nfi-1$ \KwBy $1$}
		{
    		$(\qcur_{n+1}, \pcur_{n+1}) := \C_\dt(\qcur_n, \pcur_n) + J_n$\;
			Update the relative error $\delta = E\left(\qprev,\qcur, \Nin, n+1\right)$\;
			\If{$\delta > \de$}
			{	
				$\Nfi := n$\;
				break\tcp*[l]{exit the for loop if condition satisfied; we also exit the while loop since $\delta$ is too large}
			}							
  		}
	}
}
\Output{$\{(\qcur_n, \pcur_n)\}_{0 \leq n \leq N}$, $\Nad$, $\{ k^i_{\rm conv} \}_{1 \leq i \leq \Nad}$}
\caption{Adaptive parareal algorithm \label{algo:adap}}
\end{algorithm}

\medskip

We now evaluate the total cost of the adaptive algorithm, assuming that the communication time is zero. Each time-slab $[\Nin^i \, \dt, \Nfi^i \, \dt]$ eventually identified by the algorithm has been determined in an iterative process. We denote by $[\Nin^i \, \dt, \Nfi^{i,j} \, \dt]$ the time-slab considered at the $j$-th iteration of that process (note that only the endpoint of the time-slab depends on $j$). Denoting $m_i$ the number of iterations required to identify a sufficiently short time-slab such that convergence of the parareal iterations can be reached, we have that $\{ \Nfi^{i,j} \}_{1 \leq j \leq m_i}$ is a decreasing sequence with $\Nfi^{i,1} = N$ and $\Nfi^{i,m_i} = \Nfi^i$. For any $1 \leq j < m_i$, the adaptive algorithm performs $k^{i,j}$ additional parareal iterations before realizing that the current time slab $[\Nin^i \, \dt, \Nfi^{i,j} \, \dt]$ is too long. For $j=m_i$, the adaptive algorithm performs $k^{i,m_i}$ additional parareal iterations before reaching convergence. The total number of iterations that have been performed to reach convergence on the $i$-th time-slab is given by $\dps k^i_{\rm conv} = \sum_{j=1}^{m_i} k^{i,j}$. The total cost of the adaptive algorithm is  
\begin{equation} \label{eq:def_cost}
\cost = \sum_{i=1}^{\Nad} \left[ (N - \Nin^i) \, C_c + \sum_{j=1}^{m_i} k^{i,j} \Big( (C_f+C_c) + (\Nfi^{i,j}-\Nin^i) \, C_c \Big) \right].
\end{equation}
The first term corresponds to the coarse propagation on the initially proposed $i$-th time slab, namely $[\Nin^i \, \dt,N \, \dt]$. The algorithm then proceeds with parareal iterations, on a slab which is possibly iteratively shortened. The wall-clock gain of the adaptive algorithm is
\begin{equation} \label{def:gain_1}
    \Gamma_{\rm adapt}(\de, \dc, N) = \frac{N \, C_f}{\cost}.
\end{equation}
%
%
If we additionally assume that the cost of the coarse propagator is negligible in comparison with the cost of the fine propagator (i.e. $N \, C_c \ll C_f$), we obtain $\Gamma_{\rm adapt}(\de, \dc, N) \approx \Gamma^{\rm ideal}_{\rm adapt}(\de, \dc, N)$, where the ideal gain is defined by
\begin{equation} \label{eq:borne_theo}
\Gamma^{\rm ideal}_{\rm adapt}(\de, \dc, N) := \frac{N}{\dsp \sum_{i=1}^{\Nad} k^i_{\rm conv}}.
\end{equation}

 
 


\subsection{Implementation in LAMMPS} \label{section:implem}

In order to provide a {\em non-intrusive} implementation of the Algorithms~\ref{algo:classical} and~\ref{algo:adap} in LAMMPS, we did not modify the source code of LAMMPS but proceeded as follows. The logic of the parareal algorithm is implemented in a so-called master code implemented in Python, while the force calculations and timestepping is carried out in LAMMPS. For each given $(q_n^k,p_n^k)$, the Python code calls LAMMPS through an API in order to request an advance of the system (using either $V_{\F}$ or $V_{\C}$) over a time range $\dt$ (by making $L$ time-steps of length $\ddt$), thereby computing $\F_\dt(q_n^k,p_n^k)$ and $\C_\dt(q_n^k,p_n^k)$. The jumps $\F_\dt(q_n^k,p_n^k) - \C_\dt(q_n^k,p_n^k)$ are then computed by the master code. In the sequential part, Python first requests LAMMPS to compute $\C_\dt(q_n^{k+1},p_n^{k+1})$ and then adds the jump to obtain $(q_{n+1}^{k+1},p_{n+1}^{k+1})$. 
This implementation does not run the fine integrations in parallel, but it already allows us to monitor the performance of the algorithm in terms of required parareal iterations.

We mentioned previously that, in order to ensure convergence of the parareal algorithm, both propagators should use the same array of random variables $\{ G_{\ell,n} \}_{0\leq \ell \leq L}$ in the scheme~\eqref{eq:scheme_debut}--\eqref{eq:scheme}, when propagating the system from time $n\dt$ to time $(n+1) \dt$. An ideal solution would be to first draw these random numbers and then to feed them to the coarse and fine propagators as required. However, in our implementation, the computation of $\F_\dt$ and $\C_\dt$ is performed within LAMMPS where there is no possibility to control the random increments used at each step of the scheme~\eqref{eq:scheme_debut}--\eqref{eq:scheme} (with $0 \leq \ell \leq L$). We can only control the seed of the random number generator. 

To circumvent this difficulty, we proceed as follows to reach the final time $T = N \, \dt$:
\begin{itemize}
    \item in Python, we draw a list of $N$ random numbers that we denote $\{ S_n \}_{1 \leq n \leq N}$ and that will be used as seeds by LAMMPS. Since LAMMPS expects the seed to be an integer number in a given range, we draw $\{ S_n \}_{1 \leq n \leq N}$ as a random sequence of i.i.d. integers uniformly distributed within that range.
    \item when we enter LAMMPS to compute $\F_\dt$ or $\C_\dt$ to integrate the system from time $(n-1) \, \dt$ to $n \, \dt$ (at any parareal iteration $k$), we provide the LAMMPS random number generator with the seed $S_n$. On the basis of that seed, the random number generator of LAMMPS provides the Gaussian increments $G_{\ell,n}$ for any $0 \leq \ell \leq L$. Because the seed is the same for $\F_\dt$ and $\C_\dt$ and at all parareal iterations, both schemes make use of the same sequence of increments $\{ G_{\ell,n} \}_{0 \leq \ell \leq L}$ (which of course remains the same at each parareal iteration).
\end{itemize}


The Python code corresponding to Algorithm~\ref{algo:adap} can be found in the GitHub repository at 
\begin{center}
{\tt https://github.com/OlgaGorynina/Parareal\_MD}
\end{center}
%

\subsection{Kinetic temperature} \label{section:eff_temp}

In our implementation of the parareal algorithm, LAMMPS is called $N$ times from Python, each of these calls asking LAMMPS to perform $L$ steps of the scheme~\eqref{eq:scheme_debut}--\eqref{eq:scheme}. This procedure generates a trajectory over a time interval of total length $N \, \dt = N \, L \, \ddt$. We explain in this section that, surprisingly enough, this procedure is not equivalent to performing directly $N \, L$ time-steps of~\eqref{eq:scheme_debut}--\eqref{eq:scheme}, and actually introduces a bias in the observed kinetic temperature. This is because the first step of the scheme (namely~\eqref{eq:scheme_debut}) reads differently from the subsequent ones (performed following~\eqref{eq:scheme}). In particular, one has to implement a carefully chosen temperature schedule in order to recover the correct kinetic temperature.

Consider the scheme which consists in making $L$ steps of the BBK scheme (that is, we start with~\eqref{eq:scheme_debut} and next perform $L-1$ steps of~\eqref{eq:scheme}). We show in Appendix~\ref{sec:kinetic_temp} (on the basis of analytical computations, see~\eqref{eq:variance_tot}, which are confirmed by the numerical results of Table~\ref{tab:temp_eff}) that the equilibrium kinetic temperature reached by the scheme in the limit $N \to \infty$ is 
\begin{equation} \label{eq:K_obs}
K_{\rm eq} = \beta^{-1} \left(1 - \dsp\cfrac{1}{2L} \right) + O(\ddt).
\end{equation}
Of course, if $L$ is very large, which is the standard regime in which LAMMPS is expected to be used, then $K_{\rm eq}$ is close to the target value $\beta^{-1}$, up to a time discretization error of the order of $O(\ddt)$. This explains why the scheme~\eqref{eq:scheme_debut}--\eqref{eq:scheme} is justified for MD simulations that are usually performed for long time horizons. But this is not our situation: we actually work with relatively small values of $L$ (in practice, in the numerical experiments described in Section~\ref{sec:num}, we even work with $L=1$). Therefore, we are in a situation where the kinetic temperature~\eqref{eq:K_obs} is quite different from $\beta^{-1}$. In order to guarantee that the scheme indeed reaches the correct equilibrium kinetic temperature, we propose to use a time-dependent temperature (which is indeed an option available in LAMMPS). In the first (resp. third) line of the time-integrator, instead of considering a fluctuating term of the form $\dsp \sqrt{2 \gamma \beta^{-1} \ddt} \, G_\ell$ (resp. $\dsp \sqrt{2 \gamma \beta^{-1} \ddt} \, G_{\ell+1}$) as in~\eqref{eq:scheme}, we use a term of the form $\dsp \sqrt{2 \gamma \beta_\ell^{-1} \ddt} \, G_\ell$ (resp. $\dsp \sqrt{2 \gamma \beta_{\ell+1}^{-1} \ddt} \, G_{\ell+1}$), where the temperature thus depends on the iterate number. More precisely, we consider the following scheme (compare with~\eqref{eq:scheme_debut}--\eqref{eq:scheme}). The first iterate of the scheme is given by 
\begin{equation} \label{eq:scheme_corrected_debut}
    \begin{cases}
        p_{1/2} = p_0 - \cfrac{\ddt}{2} \, \nabla V(q_0) - \cfrac{\ddt}{2} \, \gamma \, p_0 + \cfrac{1}{2} \, \sqrt{2 \gamma \beta_0^{-1} \ddt} \, G_0,\\
        q_1 = q_0 + \ddt \, p_{1/2}, \\
        p_1 = p_{1/2} - \cfrac{\ddt}{2} \, \nabla V(q_1) - \cfrac{\ddt}{2} \, \gamma \, p_{1/2} + \cfrac{1}{2} \, \sqrt{2 \gamma \beta_1^{-1} \ddt} \, G_1,
    \end{cases}
\end{equation}
where we recall that $(q_0, p_0)$ is the initial condition. The subsequent iterates (for $\ell \geq 1$) are given by
\begin{equation} \label{eq:scheme_corrected}
    \begin{cases}
        p_{\ell+1/2} = p_\ell - \cfrac{\ddt}{2} \, \nabla V(q_\ell) - \cfrac{\ddt}{2} \, \gamma \, p_{\ell-1/2} + \cfrac{1}{2} \, \sqrt{2 \gamma \beta_\ell^{-1} \ddt} \, G_\ell,\\
        q_{\ell+1} = q_\ell + \ddt \, p_{\ell+1/2}, \\
        p_{\ell+1} = p_{\ell+1/2} - \cfrac{\ddt}{2} \, \nabla V(q_{\ell+1}) - \cfrac{\ddt}{2} \, \gamma \, p_{\ell+1/2} + \cfrac{1}{2} \, \sqrt{2 \gamma \beta_{\ell+1}^{-1} \ddt} \, G_{\ell+1}.
    \end{cases}
\end{equation}
The first iterate~\eqref{eq:scheme_corrected_debut} differs from the next ones~\eqref{eq:scheme_corrected} in the same two ways as~\eqref{eq:scheme_debut} differs from~\eqref{eq:scheme}. We note that the temperature schedule appears in the scheme~\eqref{eq:scheme_corrected_debut}--\eqref{eq:scheme_corrected} in such a way that the fluctuating force in the last step of the iterate $\ell$ is equal to the fluctuating force in the first step of the iterate $\ell+1$, thereby preserving the implementation of the scheme as a Verlet scheme with a force including the damping and the fluctuation terms, in addition to the force field term (see Remark~\ref{rem:pourquoi_implem}). 

As shown in Appendix~\ref{sec:correction}, for each value of $L$, several schedules $\{ \beta_\ell \}_{0 \leq \ell \leq L}$ are possible in order to reach the correct equilibrium temperature. The only choice which is valid whatever $L$ is given by 
\begin{equation} \label{eq:choix_robuste}
\beta_\ell^{-1} = C_\ell \, \beta^{-1} \quad \text{with $C_0 = 3$ and $C_\ell = 1$ for any $\ell \geq 1$}.
\end{equation}
In the particular case when $L=1$ (which is the only case we consider in the numerical experiments of Section~\ref{sec:num}), another possible choice (and this is the one we make in Section~\ref{sec:num}) is
\begin{equation} \label{eq:choix_L_egal_1}
\beta_0^{-1} = \beta_1^{-1} = 2 \, \beta^{-1}.
\end{equation}
Details about the implementation in LAMMPS of such schedules are provided in Remark~\ref{rem:implem_schedule} within Appendix~\ref{sec:correction}.

\section{Numerical results} \label{sec:num}


\subsection{MD settings} \label{sec:MD_settings}

We demonstrate the efficiency of the adaptive algorithm described above by simulating the diffusion of a self-interstitial atom (SIA) in a tungsten lattice. To do so, we consider a perfect periodic lattice of tungsten atoms and insert an additional tungsten atom. This extra atom can relax to a number of equilibrium positions. Because of thermal fluctuations, the interstitial atom hops from one equilibrium state to another one in a metastable fashion (see e.g.~\cite{LelievreRoussetStoltz10,Uberuaga2018} for a comprehensive description of the MD context). We choose to work with SIA because of relatively small activation energies for diffusion, in contrast with the diffusion of other lattice defects, such as vacancies for instance, for which the activation energies are much larger. Because the activation energy is small, we can afford to run trajectories where we observe {\em several} jumps, which makes it possible to make statistical analysis on these jumps. For instance, we can compute the mean transition time with a reasonable statistical accuracy. As mentioned above, we use LAMMPS to perform these molecular dynamics simulations, and use the adaptive parareal algorithm discussed above (with $L=1$) to compute the trajectories. 

We consider a system containing 129 tungsten atoms, forming a BCC lattice (except for the interstitial atom) with periodic boundary conditions. The temperature (equal to $\beta^{-1}$ up to a multiplicative constant) is set to 2000 K and the damping parameter satisfies $\gamma^{-1} = 1$ ps. These values for temperature and damping parameter are used for all calculations in the current section. We consider two choices of time-step, $\ddt = 2$ fs and $\ddt = 0.5$ fs. 
 
\medskip 
 
To compute the forces on the atoms, we consider two types of interatomic potentials:
\begin{itemize}
\item empirical force fields based on a physically informed parameterized expression; we use the Embedded-Atom Method (EAM) potential~\cite{daw1984embedded};
\item machine-learned force fields using generic features as input to describe the chemical environment of each atom; we use Spectral Neighbor Analysis Potentials (SNAP)~\cite{thompson2015spectral}, where the parameters of the generic features are optimized using machine-learning techniques to reproduce (on some small configurations) the energies, forces, and stress tensors obtained by ab-initio computations; these potentials are denoted SNAP-6, \dots, SNAP-205, depending on the number of features. The larger the number of features, the more accurate the potential is with respect to ab-initio results, but also the more computationally expensive it is.
\end{itemize} 
Table~\ref{tabl:cost} presents typical computational times required to perform 5000 iterations of the scheme~\eqref{eq:scheme_corrected_debut}--\eqref{eq:scheme_corrected} (in a purely sequential manner), with various interatomic potentials, as measured on a laptop computer equipped with Intel(R) Core(TM) i7-10610U at 1.80GHz. We see that, on average, SNAP-205 is 175 times (resp. 2600 times) more expensive than SNAP-6 (resp. EAM).

\medskip

\begin{table}[htbp]
    \centering
    \begin{tabular}{l | c }
        \hline 
        \hline 
        Potential & Comp. time \\
        \hline 
        \hline 
        \hline 
        EAM      &  0.6923 \\
        SNAP-6   &  10.20  \\ 
        SNAP-15  &  22.65  \\
        SNAP-31  &  58.06  \\
        SNAP-56  &  151.2  \\
        SNAP-92  &  373.2  \\
        SNAP-141 &  861.2  \\ 
        SNAP-205 &  1787   \\ 
        \hline 
        \hline 
    \end{tabular}
    \caption{Computational time (in seconds) required to perform (in a sequential manner) 5000 iterations of~\eqref{eq:scheme_corrected_debut}--\eqref{eq:scheme_corrected}, with different interatomic potentials, on a standard laptop. \label{tabl:cost}}
\end{table}

\medskip

For every simulation described below, we first initialize the system with the following equilibration procedure. We consider a sample containing the 128 tungsten atoms located on a BCC lattice. A self-interstitial tungsten atom is then inserted in the sample and we minimize the energy (to drive the system in an equilibrium position). We then perform 10,000 steps of~\eqref{eq:scheme_debut}--\eqref{eq:scheme}, in a purely sequential manner and using the fine potential $V_\F$. The resulting thermalized configuration is our initial condition. Starting from there, we will next propagate the system for $N$ steps, either sequentially using $V_\F$, or in a parareal manner using the fine potential $V_\F$ and a coarse potential $V_\C$. We use the Voronoi analysis of OVITO~\cite{ovito} to identify the location of the SIA, which allows us to estimate transition times between metastable states.

The reference solution $\{(q_n^{\rm ref},p_n^{\rm ref})\}_{0 \leq n \leq N}$ is computed in a sequential manner with the SNAP-205 potential, which we denote $V^{\rm SNAP-205}_{\F}$. The corresponding reference average SIA residence time (that we compute on the basis of a single long trajectory of $N = 400,000$ time-steps) is 0.63568 ps.

We have at our disposal several potentials to be used as coarse potentials in the (adaptive) parareal algorithm. If we choose $V_{\C}$ to be close to $V^{\rm SNAP-205}_{\F}$, we may hope to need few iterations to reach convergence, but the discrepancy in term of cost between the fine and the coarse potential may be too small to observe any computational gain. From~\eqref{eq:def_cost}--\eqref{def:gain_1}, we indeed know that the gain of the parareal algorithm crucially depends on the ratio $C_f/C_c$. In what follows, we consider two strategies for the (adaptive) parareal algorithm:
\begin{itemize}
    \item Strategy I consists in using the SNAP-6 potential (denoted by $V^{\rm SNAP-6}_{\C}$) as coarse potential;
    \item Strategy II consists in using the EAM potential (denoted by $V^{\rm EAM}_{\C}$) as coarse potential.
\end{itemize}
In both strategies, the fine potential is the SNAP-205 potential $V^{\rm SNAP-205}_{\F}$. We thus have $C_f/C_c \sim 175$ in the first case and $C_f/C_c \sim 2600$ in the second case (on an Intel(R) Core(TM) i7-10610U machine).

\medskip

In the following sections, we investigate the accuracy of the parareal trajectories, first in a strong, trajectorial sense (in Section~\ref{sec:conv_traj}), second in a statistical sense (in Section~\ref{sec:conv_stat}). We conclude by discussing the observed computational gains (in Section~\ref{sec:gain}).

\subsection{Convergence of adaptive trajectories} \label{sec:conv_traj}

The objective of this section is to check the trajectorial convergence of the adaptive parareal algorithm. We consider the physical system described above and perform a trajectory consisting of $N = 1500$ time-steps (after equilibration), either using only the reference potential $V^{\rm SNAP-205}_{\F}$ or using the parareal algorithm. Along the trajectory, we register the number of time-steps (of length $\ddt)$ that the system spends in a given well before hoping to another one. The results, which have been computed with the time-step $\ddt = 2$ fs and the explosion threshold $\de = 0.35$, are presented in Table~\ref{tabl:convergence}. 

The first line corresponds to the reference trajectory. The next two lines correspond to the parareal results using the strategy I, and the last two lines correspond to the parareal results using the strategy II (for two values of the convergence threshold $\dc$).

We see that, when we use $\dc = 10^{-5}$, the parareal results (for strategy I and II) differ from the reference results. The parareal trajectory is not sufficiently close to the reference trajectory to obtain accurate results in terms of residence times (recall that the system is chaotic, so a small difference at some point of the trajectory may lead to a large difference in terms on the time spent in a metastable state). In contrast, when we set $\dc = 10^{-10}$, we observe that both parareal strategies essentially give the same SIA residence times as the reference solution. In the case I (resp. case II), the first seven (resp. first six) residence times are exactly reproduced by the parareal trajectory. With a small enough value of $\dc$, it is thus possible to obtain convergence on a long time interval (the horizon $T$ is sufficiently large to witness several exits of metastable states on $[0,T]$). 

Note that we do not observe (and actually do not expect to observe) an exact trajectorial convergence in the limit $\dc \to 0$, because of the inherent chaoticity of the trajectories: even when $\dc$ is as small as $10^{-10}$, the values of the last jumps are different between the parareal trajectories and the reference one, even though they are driven by exactly the same noise. This is often not important in practice since many commonly considered quantities of interest actually only depend on the law of the trajectories (statistical quantities) and not on the exact realization for a given noise. This is why we investigate the statistical accuracy of the parareal scheme in the next section.

\medskip

\begin{table}[htbp]
    \centering
    \begin{tabular}{c | c}
        \hline 
        \hline 
         & SIA residence times for reference solution
        \\ 
        \hline
        & [{\bf 122}, {\bf 23}, {\bf 27}, {\bf 476}, {\bf 14}, {\bf 32}, {\bf 560}, 245]
        \\  
        \hline 
        \hline 
        $\dc$ & SIA residence times for strategy I
        \\
        \hline
        $10^{-5}$		& [273, 1226]
        \\
        $10^{-10}$	& [{\bf 122}, {\bf 23}, {\bf 27}, {\bf 476}, {\bf 14}, {\bf 32}, {\bf 560}, 217, 26, 2]
        \\
        \hline 
        \hline 
        $\dc$ & SIA residence times for strategy II
        \\
        \hline 
        $10^{-5}$		&  [63, 27, 16, 36, 19, 34, 332, 972]
        \\ 
        $10^{-10}$	&  [{\bf 122}, {\bf 23}, {\bf 27}, {\bf 476}, {\bf 14}, {\bf 32}, 575, 15, 28, 31, 156]
        \\ 
        \hline 
        \hline 
    \end{tabular}
\caption{SIA residence times (in units of $\ddt$) for the reference and the parareal trajectories (for example, in the reference computations, the system spends 122 time-steps in the first well, then hops to a second well where it stays 23 time-steps, \dots). We mark in bold the residence times of the parareal trajectories that exactly agree with those of the reference trajectory ($\ddt = 2$ fs and $\de = 0.35$). \label{tabl:convergence}}
\end{table}

\subsection{Statistical analyses} \label{sec:conv_stat}


The aim of MD simulations, especially those using stochastic equations of motion, is most often to generate statistically correct trajectories, in contrast to very accurately obtaining a trajectory corresponding to a specific random number sequence. Quantities of interest include thermodynamical averages (which are computed as time averages along the trajectory) or dynamical information in a mean sense, such as average residence times in metastable states. 

In this section, we monitor the residence times of the SIA along a trajectory of $N = 100,000$ time-steps (computed with the time-step $\ddt = 2$ fs), and compute from the observed exit events an average residence time along with a confidence interval (in practice, 50 trajectories of 2,000 times-steps are computed; on the Figures~\ref{fig:stat_an_ref} to~\ref{fig:stat_an_10_-10} below, the vertical bar at $N=2,000$ corresponds to trajectories in which the SIA has never left its initial metastable state). 

We first consider the reference trajectory, computed only using the reference potential $V^{\rm SNAP-205}_{\F}$. Along the $N$ time-steps, we observe slightly more than 300 exit events. It is often the case that the residence time in the well is less than 200 time-steps, as shown on the top of Figure~\ref{fig:stat_an_ref}. On that figure, we plot the histogram of the SIA residence times (top figure) and the average of these residence times (bottom figure), as more and more transition events are taken into account to compute the average. We find that the average residence time is $T_{\rm avg} = 320.26$, with the confidence interval $[268.74;371.78]$ (here and throughout the article, confidence intervals are computed in a way such that the corresponding expectation belongs to the confidence interval with a probability of 0.95). Since we work with the time-step $\ddt = 2$ fs, this corresponds to an average residence time of 640 fs, a time consistent with the one obtained in Section~\ref{sec:MD_settings} (0.63568 ps) from a very long trajectory.

\medskip

\begin{figure}[htbp]
    \vspace{-3pt}
    \centering
    \includegraphics[width=0.7\textwidth]{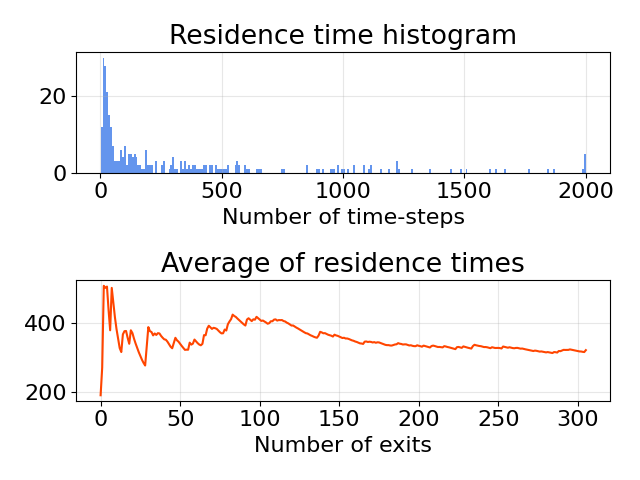}
    \caption{Reference results, in units of $\ddt$: $T_{\rm avg} = 320.26$, confidence interval $[268.74;371.78]$. 
    \label{fig:stat_an_ref}}
    \vspace{10pt}
\end{figure}

\medskip

We next repeat the same experiments, but with the EAM potential as reference potential. Results are shown on Figure~\ref{fig:stat_an_eam}. In this case, the mean residence time is $T_{\rm avg} = 99.8$ (with the confidence interval $[92.89;106.71]$). These results are very far from the reference solution results. We thus cannot rely on the EAM model to accurately predict the residence times. A coupling strategy (such as the parareal algorithm) is needed.

\begin{remark} \label{rem:stat}
Note that, in this section, all the Gaussian increments that we use (i) for the reference computations, (ii) for the EAM simulations and (iii) for the various parareal simulations presented below are independent one from each other. Likewise, the initial conditions for the three types of computations are independent. We have made this choice in view of our aim to monitor the {\em statistical} accuracy of the EAM or parareal computations with respect to the reference computations. 
On the other hand, within one type of computations (e.g. parareal simulations for given $V_\C$ and $\dc$), we have of course used the same Gaussian increments for the two propagators $\C_\dt$ and $\F_\dt$ and at all parareal iterations.
\end{remark}

\medskip

\begin{figure}[htbp]
    \vspace{-3pt}
    \centering
    \includegraphics[width=0.7\textwidth]{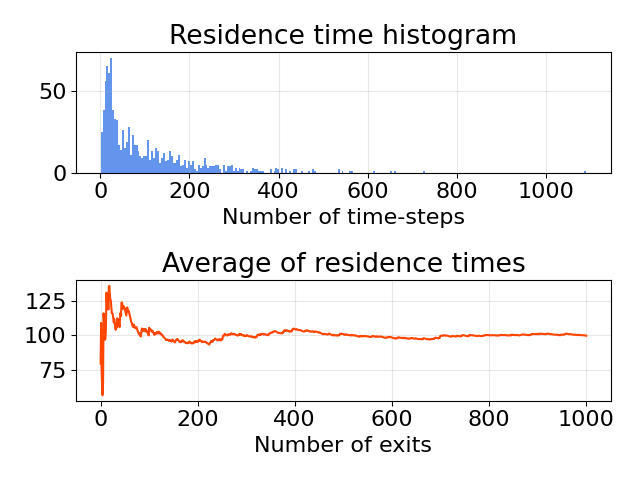}
    \caption{Results when modelling the system with the EAM potential, in units of $\ddt$: $T_{\rm avg} = 99.8$, confidence interval $[92.89;106.71]$. \label{fig:stat_an_eam}}
    \vspace{10pt}
\end{figure}

\medskip

We now consider parareal results obtained using the strategy II (coupling the fine potential $V^{\rm SNAP-205}_{\F}$ with the coarse potential $V^{\rm EAM}_{\C}$), with various values of the convergence threshold $\dc$ (for this test, we have not considered the strategy I because, as shown in Section~\ref{sec:gain} below, it provides smaller computational gains than the strategy II). Results obtained with $\dc = 10^{-3}$ (resp. $\dc = 10^{-5}$, $\dc = 10^{-10}$) are shown on Figure~\ref{fig:stat_an_10_-3} (resp. Figure~\ref{fig:stat_an_10_-5}, Figure~\ref{fig:stat_an_10_-10}). For these three values of $\dc$, we obtain confidence intervals for the mean residence time which overlap with the reference confidence interval obtained on Figure~\ref{fig:stat_an_ref}. The parareal results are thus statistically consistent with the reference results. This is true even in the case $\dc = 10^{-3}$, which is a too large value to expect trajectorial convergence. We also note that, in the case $\dc = 10^{-10}$, the statistical accuracy is {\em not} a consequence of a trajectorial convergence of the parareal trajectories to the reference trajectories. Indeed, as pointed out in Remark~\ref{rem:stat}, the initial configurations and the random noises used in the reference computations differ from those used in the parareal computations (each converged parareal trajectory is thus different from any of the reference trajectories). Moreover, for a given sequence of random noises, we observed in Section~\ref{sec:conv_traj} that the parareal trajectory differs from the reference trajectory (computed with the fine potential) after roughly 10 jumps. 

For the sake of completeness, we have also considered the agressive choice $\dc = 10^{-1}$ (results not shown). As could be expected, for this very large value of convergence threshold, the confidence intervals on residence times do not overlap at all and the parareal results are inaccurate. 

\medskip

\begin{figure}[htbp]
    \vspace{-3pt}
    \centering
    \includegraphics[width=0.7\textwidth]{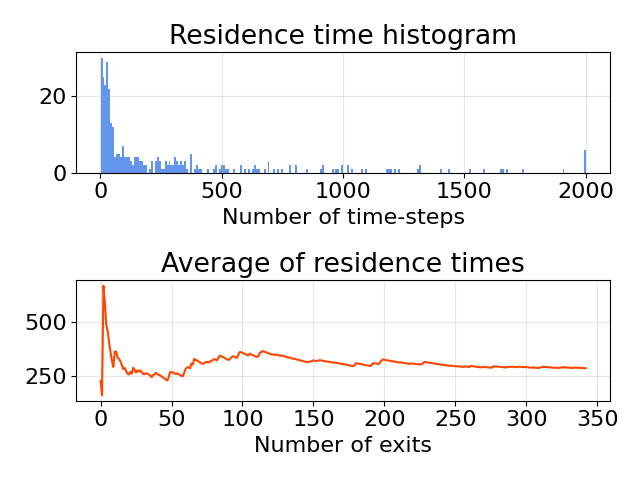}
    \caption{Parareal results (strategy II) with $\dc=10^{-3}$, in units of $\ddt$: $T_{\rm avg} = 285.714$, confidence interval $[239.73;331.7]$. \label{fig:stat_an_10_-3}}
    \vspace{10pt}
\end{figure}

\begin{figure}[htbp]
    \vspace{-3pt}
    \centering
    \includegraphics[width=0.7\textwidth]{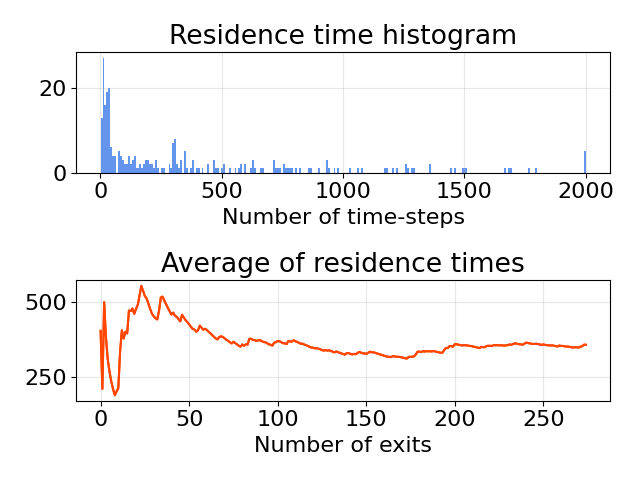}
    \caption{Parareal results (strategy II) with $\dc=10^{-5}$, in units of $\ddt$: $T_{\rm avg} = 356.36$, confidence interval $[300.9;411.83]$. \label{fig:stat_an_10_-5}}
    \vspace{10pt}
\end{figure}

\begin{figure}[htbp]
    \vspace{-3pt}
    \centering
    \includegraphics[width=0.7\textwidth]{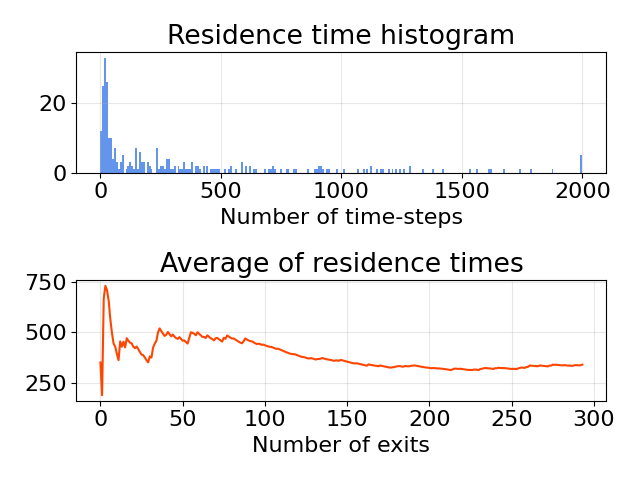}
    \caption{Parareal results (strategy II) with $\dc=10^{-10}$, in units of $\ddt$: $T_{\rm avg} = 340.14$, confidence interval $[287.82;393.45]$. \label{fig:stat_an_10_-10}}
    \vspace{10pt}
\end{figure}

\subsection{Computational gains} \label{sec:gain}

We now investigate the wall-clock gains obtained using the adaptive parareal algorithm, for various values of the time-step $\ddt$, the convergence threshold $\dc$ and the explosion threshold $\de$. We fix the time horizon at $T = N \, \ddt$ with $N = 2000$. The gain is computed using~\eqref{def:gain_1} (we also compute the ideal gain defined by~\eqref{eq:borne_theo}), where the costs are those measured on a laptop computer (see Table~\ref{tabl:cost}).


We first consider the parareal strategy I, and collect in Table~\ref{tabl:gain_snap} the values of the gain. If we use the time-step $\ddt = 2$ fs, then the maximal value of the (actual) gain is equal to 2.95 and is attained for $\de = 0.35$ and $\dc = 10^{-3}$. If we work with the smaller time-step $\ddt = 0.5$ fs, then the gain increases to 9.08 (attained for $\de = 0.3$ and $\dc = 10^{-3}$). Table~\ref{tabl:gain_eam} collects the gains for parareal strategy II. The gain reaches the value 5.18 (with the choice $\de = 0.35$ and $\dc = 10^{-3}$) when using $\ddt = 2$ fs, and increases up to 19.05 (with the choice $\de = 0.3$ and $\dc = 10^{-3}$) when using $\ddt = 0.5$ fs.

Overall, the gain obtained using strategy II is always larger than the one obtained using strategy I. This is expected since the ratio $C_f/C_c$ is more than ten times larger in the case II. In addition, the gain is always larger when considering $\ddt = 0.5$ fs rather than $\ddt = 2$ fs. 

As expected, the largest values of the gain are obtained when $\dc = 10^{-3}$ (if we decrease $\dc$ to $10^{-5}$ or $10^{-10}$, more parareal iterations are requested to achieve convergence). In terms of $\de$, the gain seems to describe a "bell shape", in the sense that it increases when $\de$ is very small and decreases when $\de$ is too large. This is in agreement with the behavior observed in the previous work~\cite{papierUpanshu}: there exists a range of values of $\de$ for which the gain remains roughly constant. 

We conclude this section by noticing that the gains obtained here are roughly similar to the gains reported in~\cite{papierUpanshu} for a Lennard-Jones cluster of 7 atoms in dimension two, although the system is here much more complex (in particular with a much larger dimensionality). 

\begin{remark}
In the case of the strategy I, we have systematically observed that the choice $\de = 0.4$ leads to a gain which is smaller than one. For some parameters choices (e.g. $\ddt = 0.5$ fs, $\de = 0.4$ and $\dc = 10^{-10}$), we decided to stop the computations when it was obvious that the gain would be smaller than one, hence the void entries in Tables~\ref{tabl:gain_snap} and~\ref{tabl:gain_eam}.

Still in the case of strategy I, we also note that the choice $\de = 0.45$ (not reported in Table~\ref{tabl:gain_snap}) may lead to unstable simulations. The system then explores regions of the phase space that are so unexpected (and so unphysical) that LAMMPS stops, declaring a {\tt lost atom}. 
\end{remark}

\begin{table}[htbp]
    \centering
    \begin{tabular}{c | l | l | c | c | c}
        \hline
        \hline
        $\ddt$ (in fs) & $\de$ & $\dc$ & $\Gamma^{\rm ideal}_{\rm adapt}$ & $\Gamma_{\rm adapt}$ & $\Nad$ \\
        \hline 
        \hline
        \hline
        0.5 & 0.15 &  $10^{-10}$	& 6.76   & 3.85 & 38 \\ 
        0.5 & 0.20 &  $10^{-10}$	& 7.75   & 4.63 & 26 \\
        0.5 & 0.25 &  $10^{-10}$	& 8.97   & 5.62 & 17 \\ 
        0.5 & 0.30 &  $10^{-10}$	& 10.27  & {\bf 8.06} & 11 \\ 
        0.5 & 0.35 &  $10^{-10}$	& {\bf 11.63}  & 3.5 & 5 \\
        0.5 & 0.40 &  $10^{-10}$	&        & <1   & \\ 
        \hline
        \hline
        0.5 & 0.15 &  $10^{-5}$	& 9.22   & 4.56 & 36 \\ 
        0.5 & 0.20 &  $10^{-5}$	& 10.7   & 5.51 & 27 \\
        0.5 & 0.25 &  $10^{-5}$	& 11.24  & {\bf 6.05} & 17 \\
        0.5 & 0.30 &  $10^{-5}$	& 10.99  & 4.44 & 9  \\ 
        0.5 & 0.35 &  $10^{-5}$	& {\bf 12.35}  & 3.77 & 4  \\ 
        0.5 & 0.40 &  $10^{-5}$	&        & <1   &    \\
        \hline
        \hline
        0.5 & 0.15 &  $10^{-3}$	& 10.53  & 4.91 & 37 \\ 
        0.5 & 0.20 &  $10^{-3}$	& 10.99  & 5.39 & 26 \\ 
        0.5 & 0.25 &  $10^{-3}$	& {\bf 13.25}  & 7.1 & 17 \\ 
        0.5 & 0.30 &  $10^{-3}$	& 13.07  & {\bf 9.08} & 9  \\ 
        0.5 & 0.35 &  $10^{-3}$	& 12.35  & 4.03 & 6  \\ 
        0.5 & 0.40 &  $10^{-3}$	&        & <1   &    \\ 
        \hline
        \hline
        2 & 0.15 &  $10^{-10}$	& 1.92  & 1.06 & 147 \\ 
        2 & 0.20 &  $10^{-10}$	& 2.24  & 1.75 & 102 \\ 
        2 & 0.25 &  $10^{-10}$	& 2.54  & 1.38 & 71  \\ 
        2 & 0.30 &  $10^{-10}$	& 2.93  & 2.26 & 42  \\ 
        2 & 0.35 &  $10^{-10}$	& {\bf 3.32}  & {\bf 2.49} & 19  \\ 
        2 & 0.40 &  $10^{-10}$	& 1.66  & 0.55 & 5   \\ 
        \hline
        \hline
        2 & 0.15 &  $10^{-5}$	& 2.84  & 1.34 & 145 \\ 
        2 & 0.20 &  $10^{-5}$	& 3.07  & 1.67 & 100 \\ 
        2 & 0.25 &  $10^{-5}$	& 3.32  & 2.1 & 68  \\ 
        2 & 0.30 &  $10^{-5}$	& 3.51  & 2.55 & 39  \\ 
        2 & 0.35 &  $10^{-5}$	& {\bf 3.67}  & {\bf 2.89} & 22  \\ 
        2 & 0.40 &  $10^{-5}$	& 2.21  & 0.81 & 7   \\ 
        \hline
        \hline
        2 & 0.15 &  $10^{-3}$	& 3.13  & 1.32 & 140 \\ 
        2 & 0.20 &  $10^{-3}$	& 3.33  & 1.76 & 102 \\ 
        2 & 0.25 &  $10^{-3}$	& 3.5   & 2.17 & 64  \\ 
        2 & 0.30 &  $10^{-3}$	& 3.54  & 2.55 & 41  \\ 
        2 & 0.35 &  $10^{-3}$	& {\bf 3.68}  & {\bf 2.95} & 24  \\ 
        2 & 0.40 &  $10^{-3}$	& 2.4   & 0.9 & 9   \\ 
        \hline
        \hline
    \end{tabular}
\caption{Wall-clock gain obtained using the adaptive parareal algorithm and strategy I (we use $V^{\rm SNAP-6}_{\C}$ and $V^{\rm SNAP-205}_{\F}$) on trajectories of length $N=2000$. We have marked in bold the best results. \label{tabl:gain_snap}}
\end{table}

\medskip

\begin{table}[htbp]
    \centering
    \begin{tabular}{c | l | l | c | c | c}
        \hline
        \hline
        $\ddt$ (in fs) & $\de$ & $\dc$ & $\Gamma^{\rm ideal}_{\rm adapt}$ &	$\Gamma_{\rm adapt}$ & $\Nad$ \\
        \hline 
        \hline
        \hline
        0.5 & 0.15 &  $10^{-10}$	& 9.9   & 9.21 & 23\\ 
        0.5 & 0.20 &  $10^{-10}$	& 12.05 & 11.48 & 14\\
        0.5 & 0.25 &  $10^{-10}$	& 13.33 & 12.88 & 9 \\ 
        0.5 & 0.30 &  $10^{-10}$	& 15.75 & {\bf 15.31} & 5 \\ 
        0.5 & 0.35 &  $10^{-10}$	& {\bf 16.39} & 14.64 & 2 \\
        0.5 & 0.40 &  $10^{-10}$	&       & <1    &   \\ 
        \hline
        \hline
        0.5 & 0.15 &  $10^{-5}$	& 14.93 & 13.99 & 22\\ 
        0.5 & 0.20 &  $10^{-5}$	& 16.26 & 15.42 & 13\\
        0.5 & 0.25 &  $10^{-5}$	& 17.24 & 16.77 & 9 \\
        0.5 & 0.30 &  $10^{-5}$	& 17.7  & 16.59 & 4 \\ 
        0.5 & 0.35 &  $10^{-5}$	& {\bf 18.52} & {\bf 17.09} & 2 \\ 
        0.5 & 0.40 &  $10^{-5}$	&       & <1    &   \\
                \hline
        \hline
        0.5 & 0.15 &  $10^{-3}$	& 16.13 & 14.88 & 22\\ 
        0.5 & 0.20 &  $10^{-3}$	& 16.95 & 15.78 & 15\\ 
        0.5 & 0.25 &  $10^{-3}$	& 17.54 & 16.81 & 9 \\ 
        0.5 & 0.30 &  $10^{-3}$	& {\bf 20.0}   & {\bf 19.05} & 5 \\ 
        0.5 & 0.35 &  $10^{-3}$	& 17.7  & 15.84 & 2 \\ 
        0.5 & 0.40 &  $10^{-3}$	&       & <1    &   \\
        \hline 
        \hline
        2 & 0.15 &  $10^{-10}$	& 2.71  & 2.58 & 89 \\ 
        2 & 0.20 &  $10^{-10}$	& 3.18  & 3.04 & 61 \\
        2 & 0.25 &  $10^{-10}$	& 3.7   & 3.56 & 40 \\
        2 & 0.30 &  $10^{-10}$	& 4.26  & 4.11 & 24 \\ 
        2 & 0.35 &  $10^{-10}$	& {\bf 4.82}  & {\bf 4.62} & 13 \\ 
        2 & 0.40 &  $10^{-10}$	& 2.27  & 1.94 & 4  \\ 
        \hline
        \hline
        2 & 0.15 &  $10^{-5}$	& 4.39  & 4.06 & 87 \\ 
        2 & 0.20 &  $10^{-5}$	& 4.56  & 4.28 & 62 \\ 
        2 & 0.25 &  $10^{-5}$	& 4.9   & 4.66 & 41 \\ 
        2 & 0.30 &  $10^{-5}$	& {\bf 5.28}  & {\bf 5.06} & 25 \\ 
        2 & 0.35 &  $10^{-5}$	& 5.22  & 4.98 & 14 \\ 
        2 & 0.40 &  $10^{-5}$	& 4.93  & 4.39 & 5  \\ 
        \hline
        \hline
        2 & 0.15 &  $10^{-3}$	& 4.73  & 4.34 & 90 \\ 
        2 & 0.20 &  $10^{-3}$	& 4.84  & 4.5 & 61 \\ 
        2 & 0.25 &  $10^{-3}$	& 5.06  & 4.8 & 38 \\ 
        2 & 0.30 &  $10^{-3}$	& 5.28  & 5.06 & 23 \\ 
        2 & 0.35 &  $10^{-3}$	& {\bf 5.45} & {\bf 5.18} & 13 \\ 
        2 & 0.40 &  $10^{-3}$	& 3.67  & 3.12 & 4  \\ 
        \hline
        \hline
    \end{tabular}
\caption{Wall-clock gain obtained using the adaptive parareal algorithm and strategy II (we use $V^{\rm EAM}_{\C}$ and $V^{\rm SNAP-205}_{\F}$) on trajectories of length $N=2000$. We have marked in bold the best results. \label{tabl:gain_eam}}
\end{table}

\appendix

\section{Kinetic temperature simulated in LAMMPS} \label{sec:app}

We consider the scheme~\eqref{eq:scheme_debut}--\eqref{eq:scheme} with $L$ time-steps, and assume that we repeat $N$ times this loop. We show in Appendix~\ref{sec:kinetic_temp} that the equilibrium kinetic temperature $K_{\rm eq}$ provided by the scheme is different from the target value $\beta^{-1}$, and that it is significantly smaller than $\beta^{-1}$ when $L$ is small. This motivates the introduction of the variant~\eqref{eq:scheme_corrected_debut}--\eqref{eq:scheme_corrected} (with the schedule~\eqref{eq:choix_robuste} or~\eqref{eq:choix_L_egal_1}), that we show in Appendix~\ref{sec:correction} to yield the correct equilibrium kinetic temperature. All the analytical derivations of this appendix are performed under the simplifying assumptions that the force field $\nabla V$ vanishes. The analytical conclusions are confirmed by numerical experiments performed on a realistic physical system modelled by a SNAP-56 potential energy.

\subsection{Stationary state of the kinetic temperature in the scheme~\eqref{eq:scheme_debut}--\eqref{eq:scheme}} \label{sec:kinetic_temp}

We consider the scheme~\eqref{eq:scheme_debut}--\eqref{eq:scheme} with a vanishing force field. We can assume without loss of generality that there is a single, one-dimensional particle: $(q,p) \in \RR \times \RR$. We define the kinetic temperature at step $\ell \in \{0,\dots,L\}$ as $K_\ell = \mathrm{Var} \, p_\ell$. We want to identify the equilibrium kinetic temperature $K_{\rm eq}$ of the scheme~\eqref{eq:scheme_debut}--\eqref{eq:scheme} with $L$ steps. This is a value such that, if $\mathrm{Var} \, p_0 = K_{\rm eq}$, then $\mathrm{Var} \, p_L = K_{\rm eq}$. This is also, by ergodicity, the limit of $\mathrm{Var} \, p_{n L}$ when $n \to \infty$.

From the first and the third lines of~\eqref{eq:scheme}, we compute that, for any $\ell \geq 1$,
\begin{align*}
    p_{\ell+1/2} 
    & = p_\ell - \cfrac{\ddt}{2} \, \gamma \, p_{\ell-1/2} + \cfrac{1}{2} \, \sqrt{2 \gamma \beta^{-1} \ddt} \, G_\ell \\
    &= \left[ p_{\ell-1/2} - \cfrac{\ddt}{2} \, \gamma \, p_{\ell-1/2} + \cfrac{1}{2} \, \sqrt{2 \gamma \beta^{-1} \ddt} \, G_\ell \right] - \cfrac{\ddt}{2} \, \gamma \, p_{\ell-1/2} + \cfrac{1}{2} \, \sqrt{2 \gamma \beta^{-1} \ddt} \, G_\ell.
\end{align*}
We set
\begin{equation} \label{eq:def_theta_mu}
    \theta = \cfrac{1}{2} \, \gamma \, \beta^{-1} \, \ddt, \qquad \mu = 1 - \cfrac{1}{2} \, \gamma \, \ddt,
\end{equation}
and we therefore have, for any $\ell \geq 2$, that
$$
p_{\ell-1/2} = (1 - \gamma \, \ddt)^{\ell-1} \, p_{1/2} + 2 \, \sqrt{\theta} \sum_{i=1}^{\ell-1} (1 - \gamma \, \ddt)^{i-1} \, G_{\ell-i}.
$$
Obviously, the above relation also holds for $\ell = 1$ (with the convention $\dps \sum_{i=1}^0 \cdot = 0$). Using the first line of~\eqref{eq:scheme_debut}, we thus deduce that, for any $\ell \geq 1$,
$$
p_{\ell-1/2} = \mu \, (1 - \gamma \, \ddt)^{\ell-1} \, p_0 + 2 \, \sqrt{\theta} \sum_{i=1}^{\ell-1} (1 - \gamma \, \ddt)^{i-1} \, G_{\ell-i} + \sqrt{\theta} \, (1 - \gamma \, \ddt)^{\ell-1} \, G_0.
$$
Using now the last line of~\eqref{eq:scheme_debut} and~\eqref{eq:scheme}, we obtain that, for any $\ell \geq 1$,
\begin{align}
     p_\ell 
     &= 
     \mu \, p_{\ell-1/2} + \sqrt{\theta} \, G_\ell 
     \nonumber
     \\
     &= \mu^2 \, (1 - \gamma \, \ddt)^{\ell-1} \, p_0 + 2 \, \mu \, \sqrt{\theta} \sum_{i=1}^{\ell-1} (1 - \gamma \, \ddt)^{i-1} \, G_{\ell-i} + \mu \, \sqrt{\theta} \, (1 - \gamma \, \ddt)^{\ell-1} \, G_0 + \sqrt{\theta} \, G_\ell.
     \label{eq:utile}
\end{align}
Using that $p_0$ and all the $G_j$, $0 \leq j \leq \ell$, are independent and centered random variables and that $\mathrm{Var} \, G_j = 1$, we compute the expectation of $p^2_\ell$ as
$$
\mathrm{Var} \, p_\ell = \mu^4 \, (1 - \gamma \, \ddt)^{2(\ell-1)} \, \mathrm{Var} \, p_0 + 4 \, \mu^2 \, \theta \sum_{i=1}^{\ell-1} (1 - \gamma \, \ddt)^{2(i-1)} + \mu^2 \, \theta \, (1 - \gamma \, \ddt)^{2(\ell-1)} + \theta.
$$
Using that $\theta = O(\ddt)$, we have, at the leading order in the time-step, and for any $\ell \geq 1$, that
\begin{align} 
    \mathrm{Var} \, p_\ell 
    &= 
    (1 - 2 \, \ell \, \gamma \, \ddt) \, \mathrm{Var} \, p_0 + 4 \, \theta \, (\ell-1) + \theta + \theta + O(\ddt^2)
    \nonumber
    \\
    &=
    (1 - 2 \, \ell \, \gamma \, \ddt) \, \mathrm{Var} \, p_0 + 4 \, \theta \left( \ell-\frac{1}{2} \right) + O(\ddt^2).
    \label{eq:variance}
\end{align}
We are now in position to identify the equilibrium kinetic temperature $K_{\rm eq}$ of the scheme~\eqref{eq:scheme_debut}--\eqref{eq:scheme} with $L$ steps. Indeed, inserting $\mathrm{Var} \, p_L = \mathrm{Var} \, p_0 = K_{\rm eq}$ in~\eqref{eq:variance}, we obtain
$$
2 \, L \, \gamma \, \ddt \, K_{\rm eq} = 4 \, \theta \left( L-\frac{1}{2} \right) + O(\ddt^2) = 2 \, \gamma \, \beta^{-1} \, \ddt \left( L-\frac{1}{2} \right) + O(\ddt^2),
$$
and therefore
\begin{equation} \label{eq:variance_tot}
    K_{\rm eq} = \beta^{-1} \left(1 - \dsp\cfrac{1}{2L} \right) + O(\ddt).
\end{equation}

We see that $K_{\rm eq}$ is always smaller than the target value $\beta^{-1}$, and that the difference is significant if $L$ is small (think again of our choice $L=1$, for which $K_{\rm eq} = \beta^{-1}/2 + O(\ddt)$). At the intermediate stages, that is for $p_\ell$ with $1 \leq \ell \leq L-1$, the result is not better: inserting~\eqref{eq:variance_tot} in~\eqref{eq:variance}, we see that, for any $\ell \in \{1,\dots,L\}$, 
\begin{equation} \label{eq:variance_middle}
    \mathrm{Var} \, p_\ell = (1 - 2 \, \ell \, \gamma \, \ddt) \, K_{\rm eq} + 4 \, \theta \left( \ell-\frac{1}{2} \right) + O(\ddt^2) = K_{\rm eq} + \gamma \, \ddt \, \beta^{-1} \left( \frac{\ell}{L} - 1 \right) + O(\ddt^2),
\end{equation}
which is close to $K_{\rm eq}$ and thus significantly different from $\beta^{-1}$ for small $L$.

\medskip 

In order to confirm the predictions of the above calculations in a more general setting (higher dimension and non-zero force field), we now turn to numerical experiments, performed using the SNAP-56 potential energy, for the same system as in Section~\ref{sec:num} (128 tungsten atoms on a BCC lattice, this time without any defect). We consider the scheme~\eqref{eq:scheme_debut}--\eqref{eq:scheme} with $L$ time-steps of length $\ddt$, and we iterate it $N$ times, in order to reach the final time $N \dt$. The numerically computed values of the equilibrium kinetic temperature $K_{\rm eq}$ are shown in Table~\ref{tab:temp_eff} for several choices of $L$, $N$ and $\beta$.

The first two lines of Table~\ref{tab:temp_eff} show that $K_{\rm eq}$ is indeed of the order of $\beta^{-1}/2$ when $L=1$, as predicted by~\eqref{eq:variance_tot}. If we set $L=10$ and $\beta^{-1} = 300$, we expect from~\eqref{eq:variance_tot} to find $K_{\rm eq}=285$. A first simulation with $N=20,000$ yields $K_{\rm eq} \approx 280$ (see third line). The discrepancy with the theoretically predicted result decreases if $N$ is increased to $N=200,000$, as shown on the fourth line (we then expect to be closer to the ergodic limit). Finally, when $L=100$ (fifth line), the difference between the computed equilibrium kinetic temperature and its target $\beta^{-1}$ is negligible.


\medskip

\begin{table}[htbp]
    \centering
    \begin{tabular}{c | c | c | c | c}
        \hline
        \hline
        $N \times L$	& $N$	& $L$	& $\beta^{-1}$	& $K_{\rm eq}$ \\
        \hline 
        \hline
        \hline
        20,000		& 20,000	& 1		    & 300			& 156.57    \\
        20,000		& 20,000	& 1		    & 600			& 303.474   \\
        200,000		& 20,000	& 10		& 300			& 280.062   \\
        2,000,000	& 200,000   & 10		& 300			& 287.56    \\
        2,000,000	& 20,000	& 100		& 300			& 303.46    \\
        \hline
        \hline
    \end{tabular}
    \caption{Equilibrium kinetic temperature obtained using the scheme~\eqref{eq:scheme_debut}--\eqref{eq:scheme}, for different values of $L$, $N$ and target temperature $\beta^{-1}$ (SNAP-56 potential energy, time-step $\ddt = 0.5$ fs, damping coefficient $\gamma^{-1} = 1$ ps). \label{tab:temp_eff}}
\end{table}

\begin{remark}
We have defined the numerical equilibrium kinetic temperature as the empirical variance of $\{ p_{n,\ell} \}_{0 \leq n \leq N, 0 < \ell \leq L}$. We could alternatively have defined it as the empirical variance of $\{ p_{n,L} \}_{0 \leq n \leq N}$. In all the test cases considered in Table~\ref{tab:temp_eff}, the first order correction in the right hand side of~\eqref{eq:variance_middle} satisfies 
$$
\left| \gamma \, \ddt \, \beta^{-1} \left( \frac{\ell}{L} - 1 \right) \right| \leq 0.3,
$$
and is thus negligible in comparison to the leading order term $\dsp K_{\rm eq} = \beta^{-1} \left( 1 - \cfrac{1}{2L} \right)$ of~\eqref{eq:variance_middle}. We thus expect the variances of $\{ p_{n,\ell} \}_{0 \leq n \leq N, 0 < \ell \leq L}$ and of $\{ p_{n,L} \}_{0 \leq n \leq N}$ to be close. The analytical result~\eqref{eq:variance_tot} corresponds to the theoretical variance of $\{ p_{n,L} \}_{0 \leq n \leq N}$.
\end{remark}

\subsection{Correction procedure} \label{sec:correction}


In order to ensure that the kinetic temperature at every time-step of our integrator remains fixed to $\beta^{-1}$, we have introduced in Section~\ref{section:eff_temp} the variant~\eqref{eq:scheme_corrected_debut}--\eqref{eq:scheme_corrected} of~\eqref{eq:scheme_debut}--\eqref{eq:scheme}. We are now going to show how to choose the temperature schedule in order to guarantee that the equilibrium kinetic temperature is indeed $K_{\rm eq} = \beta^{-1}$ at all the steps $\ell \in \{0,\dots,L\}$. 

Instead of working with $\beta_\ell$, we work with the correction constants $C_\ell$ defined by
$$
\beta_\ell^{-1} = C_\ell \, \beta^{-1}. 
$$
Similarly to~\eqref{eq:def_theta_mu}, we set
$$
\theta_\ell = \cfrac{1}{2} \, \gamma \, \beta_\ell^{-1} \, \ddt,
$$
and similarly to~\eqref{eq:utile}, we have, for any $\ell \geq 1$, that
$$
p_\ell 
= 
\mu^2 \, (1 - \gamma \, \ddt)^{\ell-1} \, p_0 + 2 \, \mu \sum_{i=1}^{\ell-1} (1 - \gamma \, \ddt)^{i-1} \, \sqrt{\theta_{\ell-i}} \, G_{\ell-i} + \mu \, \sqrt{\theta_0} \, (1 - \gamma \, \ddt)^{\ell-1} \, G_0 + \sqrt{\theta_\ell} \, G_\ell.
$$
We next compute the expectation of $p^2_\ell$ as
$$
\mathrm{Var} \, p_\ell = \mu^4 (1 - \gamma \, \ddt)^{2(\ell-1)} \, \mathrm{Var} \, p_0 + 4 \, \mu^2 \sum_{i=1}^{\ell-1} (1 - \gamma \, \ddt)^{2(i-1)} \, \theta_{\ell-i} + \mu^2 \, \theta_0 \, (1 - \gamma \, \ddt)^{2(\ell-1)} + \theta_\ell.
$$
Using that $\theta_j = O(\ddt)$ for any $j$, we have, at the leading order in the time-step, and for any $\ell \geq 1$, that
\begin{equation} 
    \mathrm{Var} \, p_\ell 
    = 
    (1 - 2 \, \ell \, \gamma \, \ddt) \, \mathrm{Var} \, p_0 + 4 \sum_{i=1}^{\ell-1} \theta_i + \theta_0 + \theta_\ell + O(\ddt^2).
    \label{eq:variance_bis}
\end{equation}
We now wish to choose $C_0$, $C_1$, \dots, $C_L$ such that, if $\mathrm{Var} \, p_0 = \beta^{-1}$, then $\mathrm{Var} \, p_\ell = \beta^{-1}$ for any $1 \leq \ell \leq L$. We thus have $L+1$ unknowns for $L$ equations. 

Setting $\ell=1$ in~\eqref{eq:variance_bis} and imposing $\mathrm{Var} \, p_0 = \mathrm{Var} \, p_1 = \beta^{-1}$ there, we get
$$
2 \, \gamma \, \ddt \, \beta^{-1} = \theta_0 + \theta_1 + O(\ddt^2) = \cfrac{1}{2} \, \gamma \, \beta^{-1} \, \ddt \, (C_0 + C_1) + O(\ddt^2),
$$
which leads to enforcing 
\begin{equation} \label{eq:conc_1}
C_1 = 4-C_0.
\end{equation}
We next infer from~\eqref{eq:variance_bis} that, for any $\ell \geq 2$,
$$
\mathrm{Var} \, p_\ell = \mathrm{Var} \, p_{\ell-1} - 2 \, \gamma \, \ddt \, \mathrm{Var} \, p_0 + 4 \, \theta_{\ell-1} + \theta_\ell - \theta_{\ell-1} + O(\ddt^2).
$$    
Imposing there that $\mathrm{Var} \, p_0 = \mathrm{Var} \, p_{\ell-1} = \mathrm{Var} \, p_\ell = \beta^{-1}$ leads to
$$
2 \, \gamma \, \ddt \, \beta^{-1} = 3 \, \theta_{\ell-1} + \theta_\ell + O(\ddt^2) = \cfrac{1}{2} \, \gamma \, \beta^{-1} \, \ddt \, (3 \, C_{\ell-1} + C_\ell) + O(\ddt^2),
$$
which leads to enforcing 
\begin{equation} \label{eq:conc_2}
C_\ell = 4 - 3 \, C_{\ell-1} \quad \text{for any $\ell \geq 2$}.
\end{equation}
Collecting~\eqref{eq:conc_1} and~\eqref{eq:conc_2}, we obtain
$$
C_\ell = 1 - (-3)^\ell - (-3)^{\ell-1} \, C_0 \quad \text{for any $1 \leq \ell \leq L$},
$$
from which we infer the simpler expression
$$
C_\ell = 1 + 3^{\ell-1} \, (C_0 - 3) \ \ \text{if $\ell$ is even}, \qquad C_\ell = 1 + 3^{\ell-1} \, (3 - C_0) \ \ \text{if $\ell$ is odd}.
$$
We now observe that not all choices of $C_0$ are admissible choices, since we have to ensure that $C_\ell > 0$ for any $0 \leq \ell \leq L$. If $\ell$ can take arbitrary large values, then the only possible choice is 
\begin{equation} \label{eq:robust}
    C_0 = 3, \qquad C_\ell = 1 \ \ \text{for any $\ell \geq 1$}.
\end{equation}
This is the only choice which is robust with respect to $L$. In contrast, if $L$ is fixed beforehand, several choices are possible. For instance, in the case $L=1$, in addition to the choice $C_0 = 3$ and $C_1 = 1$, another possible choice is $C_0 = C_1 = 2$.

\begin{remark} \label{rem:implem_schedule}
We have proceeded as follows to implement the scheme~\eqref{eq:scheme_corrected_debut}--\eqref{eq:scheme_corrected} in LAMMPS. In LAMMPS, the temperature schedule is defined in terms of the total amount of time (here, $n \, \dt + \ell \, \ddt$ for some $n$ and $\ell$) elapsed since the initial time $t=0$. In contrast, in our scheme, the effective temperature $\beta_\ell^{-1}$ depends on $\ell$ but is independent of $n$. Our implementation relies on using the command {\tt reset\_timestep} of LAMMPS, which allows to reset the simulation clock to zero and that we call whenever $t = n \, \dt$ for some $n$. Thus, LAMMPS always makes use of the temperature at time 0, namely $\beta_0^{-1}$, when considering the first line of~\eqref{eq:scheme_corrected_debut} to advance from time $n \, \dt$ to time $n \, \dt + \ddt$.
\end{remark}

\medskip

In order to check that the temperature schedules derived above indeed enforce the correct temperature in the system, we now turn to numerical experiments, which are performed with the same physical system as in Appendix~\ref{sec:kinetic_temp}. We consider the scheme~\eqref{eq:scheme_corrected_debut}--\eqref{eq:scheme_corrected} with $L$ time-steps of length $\ddt$, and we iterate it $N$ times to reach the final time $N \, \dt$. The numerically computed value of the equilibrium kinetic temperature $K_{\rm eq}$ is shown in Table~\ref{tab:temp_cor} for several choices of $L$, $N$ and correction procedure.

In the first two lines of Table~\ref{tab:temp_cor}, we set $L=1$ and consider two possible correction procedures, $C_0 = C_1 = 2$ in the first line and $C_0 = 3$, $C_1 = 1$ in the second line. We observe that both procedures lead to a numerical equilibrium kinetic temperature very close to its target (compare with the first line of Table~\ref{tab:temp_eff}). For larger values of $L$ (here, $L=10$; we have considered two values of $N$ for the sake of comparison with Table~\ref{tab:temp_eff}), we only consider the robust choice~\eqref{eq:robust}, which leads to excellent results. 


\medskip

\begin{table}[htbp]
    \centering
    \begin{tabular}{c | c | c | c | c }
        \hline
        \hline
        $N \times L$	& $N$	& $L$ & correction procedure & $K_{\rm eq}$ \\
        \hline 
        \hline
        \hline
        20,000   & 20,000    & 1  &  $C_0 = C_1 = 2$     & 303.38 \\  %
        20,000   & 20,000    & 1  &  $C_0 = 3$, $C_1 = 1$  & 296.87 \\  %
        200,000  & 20,000    & 10 &  $C_0 = 3$, $C_{[1,\dots,10]} = 1$  & 303.12 \\  
        2,000,000 & 200,000   & 10 &  $C_0 = 3$, $C_{[1,\dots,10]} = 1$  & 303.85 \\  
        \hline
        \hline
    \end{tabular}
    \caption{Equilibrium kinetic temperature obtained using the corrected scheme~\eqref{eq:scheme_corrected_debut}--\eqref{eq:scheme_corrected}, for different values of $L$ and $N$ and different correction procedures (SNAP-56 potential energy, target temperature $\beta^{-1} = 300$ K, time-step $\ddt = 0.5$ fs, damping coefficient $\gamma^{-1} = 1$ ps). \label{tab:temp_cor}}
\end{table}

\section*{Acknowledgments} Part of this work has been completed while DP was visiting Paris-Est Sup as an invited professor. The hospitality of that institution is gratefully acknowledged. This project has received funding from the Agence Nationale de la Recherche (ANR, France) and the European High-Performance Computing Joint Undertaking (JU) under grant agreement 955701 (project Time-X). The JU receives support from the European Union’s Horizon 2020 research and innovation programme and Belgium, France, Germany, Switzerland. 
This work was also partially funded by the European Research Council (ERC) under the European Union's Horizon 2020 research and innovation programme (grant 810367; project EMC2).

\bibliographystyle{plain}
\bibliography{biblio_parareal_MD}

\begin{thebibliography}{10}

\bibitem{bal2003parallelization}
G.~Bal.
\newblock Parallelization in time of (stochastic) ordinary differential
  equations.
\newblock Preprint available at {\tt
  https://www.stat.uchicago.edu/$\sim$guillaumebal/PAPERS/paralleltime.pdf}.

\bibitem{bal-maday-02}
G.~Bal and Y.~Maday.
\newblock A parareal time discretization for nonlinear {PDE}'s with application
  to the pricing of an {A}merican put.
\newblock In L.F. Pavarino and A.~Toselli, editors, {\em Recent developments in
  domain decomposition methods}, volume~23 of {\em Lecture Notes in
  Computational Science and Engineering}, pages 189--202. Springer Verlag,
  2002.

\bibitem{BBK}
A.~Blouza, L.~Boudin, and S.-M. Kaber.
\newblock Parallel in time algorithms with reduction methods for solving
  chemical kinetics.
\newblock {\em Communications in Applied Mathematics and Computational
  Science}, 5:241--263, 2010.

\bibitem{brunger1984stochastic}
A.~Br{\"u}nger, C.L. Brooks~III, and M.~Karplus.
\newblock Stochastic boundary conditions for molecular dynamics simulations of
  {ST2} water.
\newblock {\em Chemical Physics Letters}, 105(5):495--500, 1984.

\bibitem{DaiBrisLegollMaday13}
X.~Dai, C.~Le~Bris, F.~Legoll, and Y.~Maday.
\newblock Symmetric parareal algorithms for {H}amiltonian systems.
\newblock {\em Mathematical Modelling and Numerical Analysis}, 47(3):717--742,
  2013.

\bibitem{dai_maday}
X.~Dai and Y.~Maday.
\newblock Stable parareal in time method for first- and second-order hyperbolic
  systems.
\newblock {\em SIAM Journal on Scientific Computing}, 35:A52--A78, 2013.

\bibitem{daw1984embedded}
M.S. Daw and M.I. Baskes.
\newblock Embedded-atom method: Derivation and application to impurities,
  surfaces, and other defects in metals.
\newblock {\em Physical Review B}, 29(12):6443, 1984.

\bibitem{engblom2009parallel}
S.~Engblom.
\newblock Parallel in time simulation of multiscale stochastic chemical
  kinetics.
\newblock {\em SIAM Multiscale Modeling and Simulation}, 8:46--68, 2009.

\bibitem{farhat2003time}
C.~Farhat and M.~Chandesris.
\newblock Time-decomposed parallel time-integrators: theory and feasibility
  studies for fluid, structure, and fluid--structure applications.
\newblock {\em International Journal for Numerical Methods in Engineering},
  58:1397--1434, 2003.

\bibitem{gander_50years}
M.J. Gander.
\newblock 50 years of time parallel time integration.
\newblock In T.~Carraro, M.~Geiger, S.~K\"orkel, and R.~Rannacher, editors,
  {\em Multiple Shooting and Time Domain Decomposition}, pages 69--114.
  Springer Verlag, 2015.

\bibitem{gander-lunet}
M.J. Gander, T.~Lunet, D.~Ruprecht, and R.~Speck.
\newblock A unified analysis framework for iterative parallel-in-time
  algorithms.
\newblock arXiv preprint 2203.16069.

\bibitem{gander2007analysis}
M.J. Gander and S.~Vandewalle.
\newblock Analysis of the parareal time-parallel time-integration method.
\newblock {\em SIAM Journal on Scientific Computing}, 29:556--578, 2007.

\bibitem{garrido2005}
I.~Garrido, M.~Espedal, and G.~Fladmark.
\newblock A convergent algorithm for time parallelization applied to reservoir
  simulation.
\newblock In R.~Kornhuber, R.~Hoppe, J.~P\'eriaux, O.~Pironneau, O.~Widlund,
  and J.~Xu, editors, {\em Domain decomposition methods in science and
  engineering}, volume~40 of {\em Lecture Notes in Computational Science and
  Engineering}, pages 469--476. Springer Berlin Heidelberg, 2005.

\bibitem{LegollLelievreSamaey20}
F.~Legoll, T.~Leli{\`e}vre, K.~Myerscough, and G.~Samaey.
\newblock Parareal computation of stochastic differential equations with
  time-scale separation: a numerical convergence study.
\newblock {\em Computing and Visualization in Science}, 23(9), 2020.

\bibitem{LegollLelievreSamaey13}
F.~Legoll, T.~Leli{\`e}vre, and G.~Samaey.
\newblock A micro-macro parareal algorithm: application to singularly perturbed
  ordinary differential equations.
\newblock {\em SIAM Journal on Scientific Computing}, 35(4):A1951--A1986, 2013.

\bibitem{papierUpanshu}
F.~Legoll, T.~Leli\`evre, and U.~Sharma.
\newblock An adaptive parareal algorithm: application to the simulation of
  molecular dynamics trajectories.
\newblock {\em SIAM Journal on Scientific Computing}, 44(1):B146--B176, 2022.

\bibitem{LelievreRoussetStoltz10}
T.~Leli\`{e}vre, M.~Rousset, and G.~Stoltz.
\newblock {\em Free Energy Computations}.
\newblock Imperial College Press, 2010.

\bibitem{lions2001resolution}
J.-L. Lions, Y.~Maday, and G.~Turinici.
\newblock R{\'e}solution d'{EDP} par un sch{\'e}ma en temps parar{\'e}el [{A}
  "parareal" in time discretization of {PDE}'s].
\newblock {\em C. R. Acad. Sci. Paris, S\'erie I}, 332(7):661--668, 2001.

\bibitem{maday41parareal}
Y.~Maday.
\newblock Parareal in time algorithm for kinetic systems based on model
  reduction.
\newblock In A.~Bandrauk, M.C. Delfour, and C.~Le~Bris, editors, {\em
  High-dimensional partial differential equations in science and engineering},
  volume~41 of {\em CRM Proceedings and Lecture Notes}, pages 183--194.
  American Mathematical Society, 2007.

\bibitem{MulaMaday20}
Y.~Maday and O.~Mula.
\newblock An adaptive parareal algorithm.
\newblock {\em Journal of Computational and Applied Mathematics}, 377:112915,
  2020.

\bibitem{sall2016}
G.~Pag\`es, O.~Pironneau, and G.~Sall.
\newblock The parareal algorithm for {A}merican options.
\newblock {\em C. R. Acad. Sci. Paris, S\'erie I}, 354(11):1132–--1138, 2016.

\bibitem{ovito}
A.~Stukowski.
\newblock Visualization and analysis of atomistic simulation data with {OVITO}
  -- the {O}pen {V}isualization {T}ool.
\newblock {\em Modeling and Simulation in Materials Science and Engineering},
  18(1):015012, 2010.

\bibitem{LAMMPS}
A.P. Thompson, H.M. Aktulga, R.~Berger, D.S. Bolintineanu, W.M. Brown, P.S.
  Crozier, P.J. in~'t Veld, A.~Kohlmeyer, S.G. Moore, T.D. Nguyen, R.~Shan,
  M.J. Stevens, J.~Tranchida, C.~Trott, and S.J. Plimpton.
\newblock {LAMMPS} -- a flexible simulation tool for particle-based materials
  modeling at the atomic, meso, and continuum scales.
\newblock {\em Computer Physics Communications}, 271:108171, 2022.

\bibitem{thompson2015spectral}
A.P. Thompson, L.P. Swiler, C.R. Trott, S.M. Foiles, and G.J. Tucker.
\newblock Spectral neighbor analysis method for automated generation of
  quantum-accurate interatomic potentials.
\newblock {\em Journal of Computational Physics}, 285:316--330, 2015.

\bibitem{Uberuaga2018}
B.P. Uberuaga and D.~Perez.
\newblock Computational methods for long-timescale atomistic simulations.
\newblock In W.~Andreoni and S.~Yip, editors, {\em Handbook of Materials
  Modeling: Method: Theory and Modeling}, pages 683--688. Springer
  International Publishing, Cham, 2020.

\bibitem{zamora2020accelerated}
R.J. Zamora, D.~Perez, E.~Martinez, B.P. Uberuaga, and A.F. Voter.
\newblock Accelerated molecular dynamics methods in a massively parallel world.
\newblock In W.~Andreoni and S.~Yip, editors, {\em Handbook of Materials
  Modeling: Methods: Theory and Modeling}, pages 745--772. Springer
  International Publishing, Cham, 2018.

\end{thebibliography}

\end{document}